\documentclass{amsart}
\usepackage{amssymb}
\usepackage{latexsym}
\usepackage{amsfonts}
\usepackage{amsmath}

\newcommand{\bbR}{{\mathbb R}}
\newcommand{\bbC}{{\mathbb C}}
\newcommand{\bbF}{{\mathbb F}}
\newcommand{\bbK}{{\mathbb K}}
\newcommand{\bbL}{{\mathbb L}}
\newcommand{\bmu}{\boldsymbol{\mu}}
\newcommand{\bnu}{\boldsymbol{\nu}}
\newcommand{\blambda}{\boldsymbol{\lambda}}
\newcommand{\bgamma}{\boldsymbol{\gamma}}
\newcommand{\tX}{\tilde{X}}
\newcommand{\tY}{\tilde{Y}}
\newcommand{\tF}{\tilde{F}}
\newcommand{\tK}{\tilde{K}}
\newcommand{\tC}{\tilde{C}}
\newcommand{\tP}{\tilde{P}}
\newcommand{\ts}{\tilde{s}}
\newcommand{\ta}{\tilde{a}}
\newcommand{\tkappa}{\tilde{\kappa}}
\newcommand{\tLambda}{\tilde{\Lambda}}
\newcommand{\ch}{\mathrm{ch}}
\newcommand{\tch}{\tilde{\mathrm{ch}}}
\newcommand{\cP}{\mathcal{P}}
\newcommand{\cX}{\mathcal{X}}
\newcommand{\GL}{{\rm GL}}
\newcommand{\Or}{{\rm O}}
\newcommand{\SO}{{\rm SO}}

\newcommand{\U}{{\rm U}}
\newcommand{\tPhi}{\tilde{\Phi}}
\newcommand{\tTheta}{\tilde{\Theta}}

\theoremstyle{plain}
\newtheorem{theorem}{Theorem}[section]
\newtheorem{corollary}{Corollary}[section]
\newtheorem{lemma}{Lemma}[section]
\newtheorem{proposition}{Proposition}[section]

\theoremstyle{definition}

\begin{document}
\title[Extending real-valued characters]
{Extending real-valued characters of finite general linear and unitary groups on elements related to regular unipotents}

\author{Rod Gow}
\address{School of Mathematical Sciences\\
University College\\
Belfield, Dublin 4\\
Ireland}
\email{rod.gow@ucd.ie}
\author{C. Ryan Vinroot}
\address{Mathematics Department\\
The University of Arizona\\
Tucson, AZ  85721-0089\\
USA}
\email{vinroot@math.arizona.edu}
\begin{abstract} Let $\GL(n, \bbF_q)\langle \tau \rangle$ and $\U(n, \bbF_{q^2})\langle \tau \rangle$ denote the finite general linear and unitary groups extended by the transpose inverse automorphism, respectively, where $q$ is a power of $p$.  Let $n$ be odd, and let $\chi$ be an irreducible character of either of these groups which is an extension of a real-valued character of $\GL(n, \bbF_q)$ or $\U(n, \bbF_{q^2})$.  Let $y\tau$ be an element of $\GL(n, \bbF_q)\langle \tau \rangle$ or $\U(n, \bbF_{q^2})\langle \tau \rangle$ such that $(y\tau)^2$ is regular unipotent in $\GL(n, \bbF_q)$ or $\U(n, \bbF_{q^2})$, respectively.  We show that $\chi(y\tau) = \pm 1$ if $\chi(1)$ is prime to $p$ and $\chi(y\tau) = 0$ otherwise.  Several intermediate results on real conjugacy classes and real-valued characters of these groups are obtained along the way.
\end{abstract}

\maketitle
\section{Introduction}
\noindent Let $\bbF$ be a field and let $n$ be a positive integer. Let  $\GL(n, \bbF)$ denote the general linear group of degree $n$ over $\bbF$. In the special case that 
$\bbF$ is the finite field of order $q$, we denote the corresponding general linear group
by $\GL(n, \bbF_q)$. Let $\tau$
denote the involutory automorphism of $\GL(n, \bbF)$ which maps an element
$g$ to its transpose inverse $(g')^{-1}$, where $g'$ denotes the transpose of $g$, and let $\GL(n, \bbF)\langle \tau\rangle$ denote the semidirect product of
$\GL(n, \bbF)$ by $\tau$. Thus in $\GL(n, \bbF)\langle \tau\rangle$, we have $\tau^2=1$ and $\tau g\tau=(g')^{-1}$ for $g\in \GL(n, \bbF)$.
Let $g\tau$ and $h\tau$ be elements in the coset $\GL(n, \bbF) \tau$. These elements are conjugate in $\GL(n, \bbF)\langle \tau\rangle$ if and only if there is an element $x$ in
$\GL(n, \bbF)$ with
\[
xg\tau x^{-1}=h\tau,
\]
which is equivalent to the equality 
\[
xgx'=h.
\]
Identifying
$g$ and $h$ with non-degenerate bilinear forms over $\bbF$, we see that $g\tau$ and $h\tau$ are conjugate precisely when $g$ and $h$ define equivalent bilinear forms. Furthermore, it is clear that the centralizer of
$g\tau$ in $\GL(n, \bbF)$ consists of those elements $z\in \GL(n, \bbF)$ which satisfy
\[
zgz'=g.
\]
This means that we may identify the centralizer of $g\tau$ with the isometry group of the bilinear form
defined by $g$. Thus the study of the conjugacy classes and their centralizers of elements in the coset
$\GL(n, \bbF)\tau$ encompasses one of the classical problems of linear algebra.

In this paper, we are interested in the irreducible characters and conjugacy classes of the group
$\GL(n, \bbF_q)\langle \tau\rangle$. The first-named author above showed that,
when $q$ is a power of an odd prime, all the complex characters of $\GL(n, \bbF_q)\langle \tau\rangle$
may be realized over the field $\bbR$ of real numbers. This property was subsequently shown to hold
for all finite fields, and it has some interesting consequences for the characters of the finite
general linear group. 

Now the finite unitary group $\U(n, \bbF_{q^2})$ also admits
the transpose inverse map as an involutory automorphism and we may thus form a 
corresponding semidirect
product $\U(n,\bbF_{q^2})\langle \tau\rangle$.  In Section \ref{conjsec}, we will show that there is a one-to-one correspondence
between the conjugacy classes of $\GL(n, \bbF_q)\langle \tau\rangle$ in the coset $\GL(n, \bbF_q)\tau$ and the conjugacy classes of $\U(n, \bbF_{q^2})\langle \tau\rangle$ in the coset $\U(n, \bbF_{q^2})\tau$, and that this correspondence preserves the order
of elements in corresponding classes and maps centralizers to isomorphic centralizers.
Some consequences of this correspondence were described in the paper of the second named author
above and N. Thiem \cite{TV}.

In Section \ref{characters}, we give a combinatorial description of the irreducible characters of $\GL(n, \bbF_q)$ and $\U(n, \bbF_{q^2})$, and apply this description to give a correspondence between real-valued irreducible characters of these two groups.  In Section \ref{charduality}, this correspondence of characters is applied, along with the duality between semisimple and regular characters, to count the number of real-valued semisimple and regular characters in these groups.

In Section \ref{realpropsec}, we return to the subject of conjugacy.  In particular, results are obtained on which elements in $\U(n, \bbF_{q^2})$ and $\U(n, \bbF_{q^2})\langle \tau \rangle$ are strongly real.  In Lemma \ref{unipotentsquare}, we prove that there are elements $x\tau$ and $y\tau$ in $\GL(n, \bbF_q)\tau$ and $\U(n, \bbF_{q^2})\tau$ ($q$ odd), respectively, such that, when $n$ is odd, these elements are strongly real and square to regular unipotent elements, and when $n$ is even, these elements square to the negative of a regular unipotent.

The main results on character values are proven in Sections \ref{charvaluessec} and \ref{chartwosec}, which are as follows.
\\
\\
\noindent (Theorems \ref{extensionvalues} and \ref{valueschar2}) Let $n$ be odd, $q$ a power of $p$, and let $y\tau$ be an element of $\GL(n, \bbF_q)\langle \tau \rangle$ (or $\U(n, \bbF_{q^2})\langle \tau \rangle$) such that $(y\tau)^2$ is regular unipotent, and let $\chi$ be a character of $\GL(n, \bbF_q)\langle \tau \rangle$ (or $\U(n, \bbF_{q^2})\langle \tau \rangle$) which is an extension of a real-valued character of $\GL(n, \bbF_q)$ (or $\U(n, \bbF_{q^2})$).  Then $\chi(y\tau) = \pm 1$ if $\chi(1)$ is prime to $p$, and $\chi(y\tau) = 0$ otherwise.  Also, $\chi(\tau) \equiv \pm \chi(y\tau) \pmod{p}$.
\\
\\
\indent  We also show (in Theorem \ref{complexextensionvalues}) that there is no parallel result when $n$ is even, for elements which square to $-u$, where $u$ is regular unipotent.  A key tool that is used to prove these results is a similar result due to Green, Lehrer, and Lusztig \cite{GLL} which gives the values of characters of finite groups of Lie type on regular unipotent elements.  The proof of our main result when $q$ is even, given in Section \ref{chartwosec}, uses the theory of Gelfand-Graev characters in disconnected groups due to K. Sorlin \cite{sorlin1, sorlin2}.

Feit \cite{feit1} has computed the values of cuspidal characters of $\GL(n, \bbF_q)$ extended to $\GL(n, \bbF_q)\langle \tau \rangle$, with motivation coming from the fact that characters of $\GL(n, \bbF_{q^2})$ extended by the standard Frobenius map play a key role in Shintani descent \cite{shin}.  Shintani descent is relevant in this paper as well, as it gives a correspondence between real-valued characters of $\GL(n, \bbF_q)$ and $\U(n, \bbF_{q^2})$ (perhaps the same as our correspondence in Theorem \ref{realcorr}), as studied in a more general context by Digne \cite{digne}.  Evidence suggests that Shintani descent dictates a specific relationship between the character values of $\GL(n, \bbF_q)\langle \tau \rangle$ and $\U(n, \bbF_{q^2})\langle \tau \rangle$, and this paper might be viewed as a study of this relationship on a specific type of conjugacy class.  In particular, we conjecture that if $\chi$ is a real-valued irreducible character of $\U(n, \bbF_{q^2})$ with Frobenius-Schur indicator $1$, then the values of the extensions of $\chi$ to $\U(n, \bbF_{q^2})\langle \tau \rangle$ are Galois conjugates of the values of an irreducible character of $\GL(n, \bbF_{q})\langle \tau \rangle$ extended from a real-valued irreducible character of $\GL(n, \bbF_{q})$ which is related to $\chi$ through Shintani descent.  Furthermore, we conjecture that if $\chi$ is a real-valued irreducible character of $\U(n, \bbF_{q^2})$ with Frobenius-Schur indicator $-1$, the nonzero values of the extensions of $\chi$ are purely imaginary complex numbers (see Lemma \ref{extlemma} and Corollary \ref{purelyimaginary}) which are polynomials in $\sqrt{-q}$, and the real-valued character of $\GL(n, \bbF_{q})$ which is related to $\chi$ through Shintani descent will have an extension with nonzero values obtained from those of the extension of $\chi$ by substituting $\sqrt{q}$ for $\sqrt{-q}$.  A specific case of the second part of this conjecture is mentioned at the end of Section \ref{charvaluessec}.  We note that for the simple versions of these groups and their extensions, the correspondence between character values as described above can be observed in the examples given in the Atlas of finite groups \cite{At}.

Finally, there is a conjecture of Malle \cite[p. 85]{malle} which describes a relationship between the generalized Deligne-Lusztig characters of pairs of disconnected groups which are related by twisting a Frobenius automorphism by a commuting automorphism, which includes $\GL(n, \bbF_q)\langle \tau \rangle$ and $\U(n, \bbF_{q^2})\langle \tau \rangle$.  The specific correspondence of character values as we have conjectured above would certainly give insight into Malle's conjecture in this case.
\\
\\
\noindent
{\bf Acknowledgements.  } The second author thanks the University College in Dublin for hosting a visit to the School of Mathematical Sciences in the summer of 2006, when this work began, and also K. Sorlin for email correspondence regarding the material in Section \ref{chartwosec}.

\section{A correspondence of conjugacy classes for certain finite groups of Lie type} \label{conjsec}

\noindent The aim of this section is to place the finite groups  
$\GL(n, \bbF_q)\langle \tau\rangle $ and $\U(n, \bbF_{q^2})\langle \tau\rangle$ into a more
general context where we can use the theory of algebraic groups to draw some conclusions
which apply not only to these two groups but also to a number of other important finite groups of Lie type.
We assume in this section that  $\bbK = \bar{\bbF}_q$ is a fixed algebraic closure of the finite field with $q$ elements, where $q$ is a power of the prime $p$.
Furthermore, $G$ will denote a
connected linear algebraic group over  $\bbK$.
Let $F:G\to G$ denote a (standard) Frobenius map of $G$
and let $G^F$ denote the finite subgroup of fixed points
of $F$ in $G$. Suppose that $G$ has an involutory automorphism
$\tau$ which commutes with $F$
in its action on $G$. 
We may then form a twisted
Frobenius map $\tF: G\to G$ by setting 
\[
\tF(g)=F(\tau(g))
\]
for all $g$ in $G$. Since $\tF^2=F^2$, it follows that
the subgroup $G^{\tF}$ is contained in $G^{F^2}$.

Let $x$ be any element of $G^F$. The Lang--Steinberg
theorem \cite[Theorem 10.1]{St} shows that there exists an element $z$ in $G$
with $x=z^{-1}\tF(z)$. Since $F(x)=x$, and
$F$ commutes with $\tF$, it follows that
\[
F(z)^{-1}\tF(F(z))=z^{-1}\tF(z)
\]
and thus
\[
zF(z)^{-1}=\tF(z)\tF(F(z)^{-1})=\tF(zF(z)^{-1}).
\]
This shows that the element $y=zF(z)^{-1}$ is in $G^{\tF}$.
A different choice of $z$ in $G$ used to represent $x$
according to the Lang--Steinberg theorem will lead to another
element in $G^{\tF}$ which is not obviously related to the element
$y$ just obtained. The purpose of this section is to show that
the idea of associating $x$ in $G^F$ with $y$ in $G^{\tF}$
can be used to define a correspondence of certain conjugacy
classes in extension groups of $G^F$ and $G^{\tF}$, 
respectively, as we will now explain. 

Following the construction described in the Introduction, let $G\langle \tau\rangle$ denote the semidirect product 
of $G$ by $\tau$. Since $F$ commutes with $\tau$, it follows that
both $G^F$ and $G^{\tF}$ admit $\tau$ as an automorphism
and are thus normalized by $\tau$ in $G\langle \tau\rangle$. Let
$G^F\langle \tau\rangle$ and $G^{\tF}\langle \tau\rangle$ 
denote the corresponding
subgroups of $G\langle \tau\rangle$ generated by 
$\tau$ and $G^F$, $G^{\tF}$ respectively. 

Let $H$ denote any of the groups $G$, $G^F$ or $G^{\tF}$.
In order to define our correspondence of conjugacy classes,
we prove some elementary results relating to conjugacy
of elements in $H\langle \tau\rangle$. We make use of the observation in the Introduction that
elements $a\tau$ and $b\tau$ in $H\langle \tau\rangle$ are conjugate if and only
there exists an element $c$ in $H$ with $c^{-1}a\tau(c)=b$.

\begin{lemma} \label{conjugate} Let $x$ be an element of $G^F$ and write
$x=z^{-1}\tF(z)$ for some $z\in G$. Suppose that $x\tau$ is conjugate
in $G^F\langle \tau\rangle$ to $w\tau$, where $w=v^{-1}\tF(v)$ for some $v\in G$.
Then
\[
zF(z)^{-1}\tau\ \mbox{ and }\ vF(v)^{-1}\tau
\]
are conjugate in $G^{\tF}\langle \tau\rangle$. Thus, if we also have
$x=z_1^{-1}\tF(z_1)$ for some other element $z_1$ in $G$, the elements
\[
zF(z)^{-1}\tau\ \mbox{ and }\ z_1F(z_1)^{-1}\tau
\]
are conjugate in $G^{\tF}\langle \tau\rangle$.

\end{lemma}

\begin{proof} As we noted above, there exists $g\in G^F$ with
$g^{-1}x\tau(g)=w$. Moreover, as $\tau$ is involutory and $g\in G^F$,
we have $\tF(g)=\tau(g)$. It follows that
\[
g^{-1}z^{-1}\tF(z)\tau(g)=(zg)^{-1}\tF(zg)=w=v^{-1}\tF(v).
\]
We deduce that $zgv^{-1}\in G^{\tF}$. We set
$zgv^{-1}=u$, and then obtain $v=u^{-1}zg$, where $u\in G^{\tF}$.
Since $g=F(g)$ and $\tau(u)=F(u)$, it follows that
\[
vF(v)^{-1}=u^{-1}zgF(g)^{-1}F(z)^{-1}F(u)=u^{-1}zF(z)^{-1}\tau(u),
\]
and this equality proves that $vF(v)^{-1}\tau$ and 
$zF(z)^{-1}\tau$ are conjugate in $G^{\tF}\langle \tau\rangle$, as required. The second
part is clear by taking $w=x$ and $z_1=v$.
\end{proof}

Given an element $x$ of $H$, we let $[x\tau]$ denote the conjugacy
class of $x\tau$ in $H\langle \tau\rangle$. We trust that context will make it clear
which subgroup $H$ is implied in the event of possible ambiguity.
We now define a map $\phi$ associating a conjugacy class
$[x\tau]$ in $G^F \tau$ to a conjugacy class
$[y\tau]$ in $G^{\tF} \tau$ in the following way. Write $x$ as
$z^{-1}\tF(z)$ and let $y=zF(z)^{-1}\in G^{\tF}$. Then we set
\[
\phi[x\tau]=[y\tau].
\]
Lemma \ref{conjugate} shows that the definition of $\phi$ does not depend
on the choice of $z$ to represent $x$ or the choice of $x$
to represent the conjugacy class $[x\tau]$. We note however
that $\phi$ is only defined at the level of conjugacy
classes and does not apply to individual elements.  We also note that Lemma \ref{correspondence} is a special case of \cite[Prop. 5.7]{dignemichel}.

\begin{lemma} \label{correspondence} The map $\phi$ defines 
a one-to-one correspondence between the conjugacy classes in $G^F\tau$ and the conjugacy classes
in $G^{\tF}\tau$.

\end{lemma}

\begin{proof} We first show that $\phi$ is injective.
Suppose then that $\phi[x\tau]=\phi[x_1\tau]$. Write
\[
x=z^{-1}\tF(z),\quad x_1=z_1^{-1}\tF(z_1)
\]
where $z$ and $z_1$ are appropriate elements of $G$. Then
there exists some $u\in G^{\tF}$ with
\[
u^{-1}zF(z)^{-1}\tau(u)=z_1F(z_1)^{-1}.
\]
Since $\tau(u)=F(u)$, this implies that $z^{-1}uz_1$ is in $G^F$.
We set $g=z^{-1}uz_1$. Then, since $\tau(g)=\tF(g)$,
we have
\[
g^{-1}z^{-1}\tF(z)\tau(g)=(zg)^{-1}\tF(zg)=(uz_1)^{-1}\tF(uz_1)=
z_1^{-1}\tF(z_1)
\]
and this implies that $g$ conjugates $x\tau$ into $x_1\tau$. Thus
$[x\tau]=[x_1\tau]$ and it follows that $\phi$ is injective.

Next, we show that $\phi$ is surjective. Let $u$ be any element
of $G^{\tF}$. The Lang--Steinberg theorem implies that
$u=zF(z)^{-1}$ for some $z\in G$. Since $\tF(u)=u$, we readily
check that $z^{-1}\tF(z)$ is in $G^F$. Thus if we put
$x=z^{-1}\tF(z)$, we have $\phi[x\tau]=[u\tau]$, which implies that
$\phi$ is surjective, as required.
\end{proof}

We now show that if the order of an element in $[x\tau]$ is $r$,
the order of an element in $\phi[x\tau]$ is also $r$.
Thus $\phi$ preserves the order of the elements
in a conjugacy class.

\begin{lemma} \label{conjugacy} Given $x\in G^F$, let
$[y\tau]=\phi[x\tau]$. Then $(x\tau)^{-2}$ and $(y\tau)^2$ are conjugate in $G$. Hence, 
$x\tau$ and $y\tau$ have the same
(finite) multiplicative order in  $G\langle \tau\rangle$.
\end{lemma}

\begin{proof}  As usual, we write $x=z^{-1}\tF(z)$ and set $y=zF(z)^{-1}$.
Then we have
\begin{align*}
(x\tau)^2&= x\tau(x)=z^{-1}\tF(z)(\tau(z))^{-1}F(z)\\
(y\tau)^2&= y\tau(y)=zF(z)^{-1}\tau(z) \tF(z)^{-1}.
\end{align*}
It follows that
\[
z^{-1}(y\tau)^2z=(x\tau)^{-2},
\]
as required. Furthermore,
since $x\tau$ and $y\tau$ have finite even order, and their squares have the same order by the argument above, we deduce that $x\tau$ and $y\tau$ have the same order in $G\langle \tau\rangle$.
\end{proof}

We note that the fact that $(x\tau)^{-2}$ and $(y\tau)^2$ are conjugate in $G$ (and possibly in the smaller group $G^{F^2}$) provides information on how to recognize the class $\phi[x\tau]$ in terms
of the class $[x\tau]$.

The map $\phi$ has an additional useful
property, since the centralizer of
$x\tau$ in $G^F$ is conjugate in $G$ to the centralizer
of $y\tau$ in $G^{\tF}$, where $y\tau\in\phi[x\tau]$, as we now show.

\begin{lemma} \label{centralizers} Let $x$ be an element of $G^F$, 
with $x=z^{-1}\tF(z)$ for some $z\in G$. Let
$y=zF(z)^{-1}$. Then the centralizer of $x\tau$ in $G^F$ is
$z^{-1}G^{\tF}z\cap G^F$ and the centralizer of $y\tau$ in $G^{\tF}$ is
$zG^Fz^{-1}\cap G^{\tF}$. Thus, since these are conjugate subgroups,
the centralizer of $x\tau$ in $G^F$ is
isomorphic to the centralizer of $y\tau$ in $G^{\tF}$.
\end{lemma}

\begin{proof} An element $u\in G^F$ commutes with $x\tau$ if and only if
$u^{-1}x\tau u=x\tau$. This occurs if and only if $u^{-1}x\tau(u)=x$.
Since $\tau(u)=\tF(u)$, $u$ commutes with $x\tau$ if and only if
$zuz^{-1}$ is in $G^{\tF}$. Thus the centralizer of
$x\tau$ in $G^F$ is
$z^{-1}G^{\tF}z\cap G^F$, and a similar argument shows that
the centralizer of
$y\tau$ in $G^{\tF}$ is
$zG^Fz^{-1}\cap G^{\tF}$. Since
\[
z(z^{-1}G^{\tF}z\cap G^F)z^{-1}=
zG^Fz^{-1}\cap G^{\tF},
\]
the two centralizers are conjugate in $G$ and hence isomorphic.
\end{proof}

We sum up our findings related to $\phi$ in the following
theorem, which amalgamates the various lemmas we have proved.

\begin{theorem} \label{classcorrespondence} Let $G$ be a connected 
linear algebraic group over the
algebraic closure of a finite field.
Let $F:G\to G$ denote a standard Frobenius map of $G$.
Suppose that $G$ has an involutory automorphism
$\tau$ which commutes with $F$ and let $\tF$ denote
the corresponding twisted Frobenius map. Let
$H$ denote either $G^F$ or $G^{\tF}$ and 
let $H\langle \tau\rangle$ denote the semidirect product of $H$ by $\tau$.
Given $h\in H$, let $[h\tau]$ denote the conjugacy class
of $h\tau$ in $H\langle \tau\rangle$.
Given $x\in G^F$, write $x=z^{-1}\tF(z)$ for some $z\in G$ 
and set $y=zF(z)^{-1}\in G^{\tF}$.

Then the map $\phi$ defined by $\phi[x\tau]=[y\tau]$
is a one-to-one correspondence between the conjugacy classes
in the coset $G^F\tau$ and the conjugacy classes in the coset
$G^{\tF} \tau$. 
The elements $x\tau$ and $y\tau$ have the same
order and the centralizer of $x\tau$ in $G^F$ is isomorphic
to the centralizer of $y\tau$ in $G^{\tF}$.
\end{theorem}

We note that our theorem applies when we take $G$ to be the group $\GL(n, \bbK)$, where
$\bbK$ is the algebraic closure of a finite field and $\tau$ is the transpose inverse automorphism.
The corresponding groups $G^F$ and  $G^{\tF}$ are $\GL(n, \bbF_q) $ and $\U(n, \bbF_{q^2})$, which will be the main application in this paper.  In this case, Theorem \ref{classcorrespondence} says that the number of conjugacy classes in $\GL(n, \bbF_q)\tau$ is equal to the number of conjugacy classes in $\U(n, \bbF_{q^2})\tau$, and the centralizers in $\GL(n, \bbF_q)$ and $\U(n, \bbF_{q^2})$ of corresponding classes  are isomorphic.  There are several results for the number of conjugacy classes in $\GL(n, \bbF_q)\tau$ and the sizes of their $\GL(n, \bbF_q)$-centralizers given by Fulman and Guralnick in \cite[Sections 6 and 9]{fulman}.  By applying Theorem \ref{classcorrespondence}, we obtain identical results for conjugacy classes in $\U(n, \bbF_{q^2})\tau$.

We mention another special case of our theorem,
which does not appear to be obvious from the standpoint
of linear algebra.  We take $G$ to be the special orthogonal group $\SO(2m, \bbK)$ of even degree
$2m$ over $\bbK$. $G$ is a connected linear algebraic group which has index 2 in the (disconnected)
full orthogonal group $\Or(2m,\bbK)$.  $\Or(2m,\bbK)$ contains an orthogonal reflection, $t$, say, which is an
element of order 2 not in $\SO(2m,\bbK)$. We may assume that $t$ has coefficients in the field
of order $p$. Conjugation by $t$ induces an involutory automorphism $\tau$, say, of $\SO(2m,\bbK)$
which commutes with the Frobenius map $F$, since $t$ has coefficients in the prime field. We may thus also form a twisted Frobenius map $\tF$ by means of $\tau$. The finite group $G^F$ is then the split special orthogonal group $\SO^+(2m,\bbF_q)$ and  $G^{\tF}$ is  the non-split special orthogonal group $\SO^-(2m,\bbF_q)$. We may identify the extended groups $G^F\langle \tau\rangle$ 
and $G^{\tF}\langle \tau\rangle$ with the full orthogonal groups $\Or^+(2m,\bbF_q)$ and $\Or^-(2m,\bbF_q)$, respectively. The following result summarizes how Theorem \ref{classcorrespondence} applies in this case.

\begin{corollary} There is a one-to-one correspondence
between the conjugacy classes of 
$\Or^+(2m,\bbF_q)\setminus  \SO^+(2m,\bbF_q)$ and
those of 
$\Or^-(2m,\bbF_q)\setminus \SO^-(2m,\bbF_q)$
which preserves the order of the elements in corresponding conjugacy
classes. Under this correspondence, the centralizer in 
$\SO^+(2m,\bbF_q)$
of an element in a conjugacy class in
$\Or^+(2m,\bbF_q)\setminus \SO^+(2m,\bbF_q)$
is isomorphic to the centralizer in $\SO^-(2m,\bbF_q)$
of an element in the corresponding conjugacy class of
$\Or^-(2m,\bbF_q)\setminus \SO^-(2m,\bbF_q)$.
\end{corollary}

\section{Characters of $\GL(n, \bbF_q)$ and $\U(n, \bbF_{q^2})$} \label{characters}

In this section, we give a combinatorial description of the irreducible characters of the finite general linear and unitary groups.  The development will largely follow \cite[Chapter IV]{Mac} in the case of ${\rm GL}(n, \bbF_q)$, and \cite{TV} in the case of ${\rm U}(n, \bbF_{q^2})$, where notation will vary slightly due to the fact that we give the description of characters for both cases simultaneously.

As before, we let $\bbK = \bar{\bbF}_q$ denote a fixed algebraic closure of the finite field with $q$ elements.  We set $\bar{G}_n = \GL(n, \bbK)$  and let $F: \bar{G}_n \rightarrow \bar{G}_n$ denote the standard Frobenius map defined by $F((a_{ij})) = (a_{ij}^q)$. We also let $\tF: \bar{G}_n \rightarrow \bar{G}_n$ denote the twisted Frobenius map defined by $\tF(g) = (F(g)')^{-1}$.  Then we have 
$$
 \bar{G}_n^F =  \GL(n, \bbF_q) \; \text{ and } \; \bar{G}_n^{\tF} = \U(n, \bbF_{q^2}).
$$
We also use the notation $G_n$ and $U_n$ for the groups $\GL(n, \bbF_q)$ and $\U(n, \bbF_{q^2})$, respectively.

Both $F$ and $\tF$ act on $\bar{G}_1 = \bbK^{\times}$ and the group of complex characters of $\bbK^{\times}$, which we will denote by $\hat{\bbK}^{\times}$.  We consider the orbits arising from these actions.  Let
$$
 \Phi = \{ F\text{-orbits of } \bbK^{\times} \}, \quad \tPhi = \{ \tF\text{-orbits of } \bbK^{\times} \}, 
 $$
$$ \Theta = \{ F\text{-orbits of } \hat{\bbK}^{\times} \}, \quad \tTheta = \{ \tF\text{-orbits of } \hat{\bbK}^{\times} \}. $$
{\bf Remark.} We note that the elements of $\Phi$ correspond to irreducible monic polynomials with non-zero constant term over $\bbF_q$, and there is a non-canonical bijective correspondence between the orbits in $\Phi$ and $\Theta$ (and between $\tPhi$ and $\tTheta$) which preserves the sizes of orbits.  The elements of $\tPhi$ correspond to certain monic polynomials with non-zero constant term over $\bbF_{q^2}$ which were studied and characterized by Ennola \cite{enn1}.

Let $\cP$ denote the set of partitions of non-negative integers.  For $\cX = \Phi$, $\Theta$, $\tPhi$ or $\tTheta$, we define an {\em $\cX$-partition} to be a function $\blambda: \cX \rightarrow \cP$.  The size of an $\cX$-partition is defined to be
$$||\blambda|| = \sum_{x \in \cX} |x| \blambda(x),$$
where $|x|$ denotes the cardinality of the orbit $x \in \cX$, and $|\blambda(x)|$ denotes the size of the partition $\blambda(x) \in \cP$.  Now define
$$ \cP^{\cX}_n = \{ \text{$\cX$-partition } \blambda \; \mid \; ||\blambda|| = n \} \; \; \text{ and } \;\; \cP^{\cX} = \bigcup_{n = 1}^{\infty} \cP^{\cX}_n.$$

The following parameterizations of conjugacy classes follow from the theory of elementary divisors in the case of $\GL(n, \bbF_q)$ (see \cite[IV.2]{Mac}), and follow from the work of Wall \cite{wall} and Ennola \cite{enn1} in the case of $\U(n, \bbF_{q^2})$.

\begin{theorem} The conjugacy classes $K^{\bmu}$ of $G_n$ are parameterized by $\bmu \in \cP^{\Phi}_n$ and the conjugacy classes $\tK^{\bgamma}$ of $U_n$ are parameterized by $\bgamma \in \cP^{\tPhi}_n$.
\end{theorem}

Let $C_n$ and $\tC_n$ be the rings of $\bbC$-valued class functions of $G_n$ and $U_n$, respectively, and let
$$ C = \bigoplus_n C_n \text{ and } \tC = \bigoplus_n \tC_n.$$
For $\zeta_1 \in C_i$ and $\zeta_2 \in C_j$, we define $\zeta_1 \cdot \zeta_2$ to be the class function obtained by parabolic induction, so that $\zeta_1 \cdot \zeta_2$ is obtained by inflating $\zeta_1 \otimes \zeta_2$ from $G_i \times G_j$ to the corresponding parabolic subgroup, and then inducing to $G_{i+j}$.  Then we have $\zeta_1 \cdot \zeta_2 \in C_{i+j}$.  

For $\beta_1 \in \tC_i$ and $\beta_2 \in \tC_j$, we define $\beta_1 \circ \beta_2$  to be the class function obtained by Deligne-Lusztig induction, so that $\beta_1 \circ \beta_2 = R_{U_i \times U_j}^{U_{i+j}} (\beta_1 \otimes \beta_2) \in \tC_{i+j}$ (see one of \cite{dellusz, dmbook, Ca85} for a definition of Deligne-Lusztig induction).  

With these products, $C$ and $\tC$, are graded $\bbC$-algebras with respect to the products $\cdot$ and $\circ$, respectively.  The rings $C$ and $\tC$ are endowed with the natural inner product for class functions, where $C_i$ and $C_j$, respectively $\tC_i$ and $\tC_j$, are mutually orthogonal if $i \neq j$.  This is explained in detail in \cite[Chapter IV]{Mac} in the $G_n$ case and in both \cite{dignemichelu} and \cite{TV} in the $U_n$ case. 

Let $\kappa^{\bmu} \in C_n$ be the indicator class function for the conjugacy class $K^{\bmu}$ of $G_n$, where $\bmu \in \cP^{\Phi}_n$, and similarly let $\tkappa^{\bgamma} \in \tC_n$ be the indicator class function for the conjugacy class $\tK^{\bgamma}$ of $U_n$.  We let $a_{\bmu}$ and $\ta_{\bgamma}$ denote the orders of the centralizers of elements in the conjugacy classes $K^{\bmu}$ and $\tK^{\bgamma}$, respectively.  Note that the $\kappa^{\bmu}$ and $\tkappa^{\bgamma}$ are bases of $C$ and $\tC$, respectively.

We now define two rings of symmetric functions in order to describe the irreducible characters of $G_n$ and $U_n$.  We refer to \cite[Chapter I]{Mac} for basic definitions and notions in symmetric function theory.  For each $f \in \Phi$, we let $\{X^{(f)}_i  \mid  i >0 \} = X^{(f)}$ be an infinite set of indeterminates, and similarly for each $h \in \tPhi$, we have the set of independent variables $\{ \tX^{(h)}_i  \mid i >0 \} = \tX^{(h)}$.  Let $p_n(X)$ denote the $n$-th power sum symmetric function in the set of variables $X$.  Define also for each $\varphi \in \Theta$ a set of variables $\{ Y^{(\varphi)}_i \mid i > 0\} = Y^{(\varphi)}$, and for each $\vartheta \in \tTheta$ a set of variables $\{ \tY^{(\vartheta)}_i \mid i > 0\} = \tY^{(\vartheta)}$.  We relate symmetric functions in the $X^{(f)}$ variables and the $Y^{(\varphi)}$ variables, and symmetric functions in the $\tX^{(h)}$ variables and the $\tY^{(\vartheta)}$, through the following transforms:
\begin{equation}
p_n(Y^{(\varphi)}) = (-1)^{n|\varphi| - 1} \sum_{\alpha \in \bar{G}_1^{F^{n|\varphi|}}} \xi(\alpha) p_{n|\varphi|/|f_{\alpha}|}(X^{(f_{\alpha})}),
\end{equation}
where $\alpha \in f_{\alpha}$ and $\xi \in \varphi$ for $f_{\alpha} \in \Phi$ and $\varphi \in \Theta$, and 
\begin{equation} \label{utransform}
p_n(\tY^{(\vartheta)}) = (-1)^{n|\vartheta| - 1} \sum_{\alpha \in \bar{G}_1^{\tF^{n|\vartheta|}}} \xi(\alpha) p_{n|\vartheta|/|h_{\alpha}|}(\tX^{(h_{\alpha})}),
\end{equation}
where $\alpha \in h_{\alpha}$ and $\xi \in \vartheta$ for $h_{\alpha} \in \tPhi$ and $\vartheta \in \tTheta$.   If $a$ is not a positive integer, the power symmetric function $p_a$ is defined to be $0$.  We note that these equations do not depend on the choice of $\xi$ from the $F$-orbit $\varphi$ or the $\tilde{F}$-orbit $\vartheta$.

For $\mu = (\mu_1, \mu_2, \ldots ) \in \cP$, we define $n(\mu) = \sum_i (i-1)\mu_i$, and for $\bmu \in \cP^{\cX}$, we define $n(\bmu) = \sum_{x \in \cX} |x| n(\bmu(x))$.  For $\bmu \in \cP^{\Phi}$ and $\bgamma \in \cP^{\tPhi}$, let
$$P_{\bmu} = q^{-n(\bmu)} \prod_{f \in \Phi} P_{\bmu(f)} (X^{(f)}; q^{-|f|}) \; \text{ and } \; \tP_{\bgamma} = (-q)^{-n(\bgamma)} \prod_{h \in \tPhi} P_{\bgamma(h)} (\tX^{(h)}; (-q)^{-|h|}),$$
where $P_{\lambda}(X; t)$ denotes the Hall-Littlewood symmetric function.  Now let
$$\Lambda_n = \bbC\text{-span} \{ P_{\bmu} \mid \bmu \in \cP^{\Phi}_n \} \text{ and } \tLambda_n = \bbC\text{-span} \{ P_{\bgamma} \mid \bgamma \in \cP^{\tPhi}_n \}.$$
The two rings of symmetric functions defined by 
$$\Lambda = \bigoplus_n \Lambda_n \text{ and } \tLambda = \bigoplus_n \tLambda_n$$
are graded $\bbC$-algebras with graded product given by ordinary multiplication of symmetric functions.  We define hermitian inner products on $\Lambda$ and $\tLambda$, respectively, by letting
$$\langle P_{\bmu_1}, P_{\bmu_2} \rangle = \delta_{\bmu_1 \bmu_2} a_{\bmu_1}^{-1} \; \text{ and } \; \langle \tP_{\bgamma_1}, \tP_{\bgamma_2} \rangle = \delta_{\bgamma_1 \bgamma_2} \ta_{\bgamma_1}^{-1},$$
and extending.
For $\blambda \in \cP^{\Theta}$ and $\bnu \in \cP^{\tTheta}$, we define
$$ s_{\blambda} = \prod_{\varphi \in \Theta} s_{\blambda(\varphi)}(Y^{(\varphi)}) \; \text{ and } \; \ts_{\blambda} = \prod_{\vartheta \in \tTheta} s_{\blambda(\vartheta)}(\tY^{(\vartheta)}),$$
where $s_{\lambda}(Y)$ denotes the Schur symmetric function in the set of variables $Y$.

The following theorem, in the case of the general linear groups, is due to Green \cite{green}.  In the case of the unitary groups, this theorem was originally a conjecture of Ennola \cite{enn2}, and after progress of Hotta, Springer, Lusztig, and Srinivasan \cite{hotspr, luszsrin}, it was finally proved in full generality by Kawanaka \cite{kawan}.

\begin{theorem} \label{charmap}
Define maps $\ch: C \rightarrow \Lambda$ and $\tch: \tC \rightarrow \tLambda$ by letting $\ch(\kappa_{\bmu}) = P_{\bmu}$ and $\tch(\tkappa_{\bgamma}) = \tP_{\bgamma}$, and extend linearly.  Then $\ch$ and $\tch$ are both isometric isomorphisms of graded $\bbC$-algebras.

For $\blambda \in \cP^{\Theta}$ and $\bnu \in \cP^{\tTheta}$, define the class functions
$$\chi^{\blambda} = \ch^{-1}(s_{\blambda}) \; \text{ and } \; \psi^{\bnu} = \tch^{-1}( (-1)^{\lfloor ||\bnu||/2 \rfloor + n(\bnu)} \ts_{\bnu}).$$  
Then
$$ \{ \chi^{\blambda} \; \mid \; \blambda \in \cP^{\Theta}_n \} \; \text{ and } \; \{ \psi^{\bnu} \; \mid \; \bnu \in \cP^{\tTheta}_n \} $$
are the sets of irreducible complex characters of $G_n$ and $U_n$, respectively.
\end{theorem}

Let $\chi^{\blambda}(\bmu)$ and $\psi^{\bnu}(\bgamma)$ denote the values of the characters $\chi^{\blambda}$ and $\psi^{\bnu}$ on the conjugacy classes $K^{\bmu}$ and $\tK^{\bgamma}$, respectively.  Then Theorem \ref{charmap} implies that we have
 \begin{equation} \label{expand}
 (-1)^{\lfloor ||\bnu||/2 \rfloor + n(\bnu)} \ts_{\bnu} = \sum_{\bgamma \in \cP^{\tPhi}} \psi^{\bnu}(\bgamma) \tP_{\bgamma},
 \end{equation}
 and a similar expansion for the characters of $G_n$.
 
 Let $\varphi \in \Theta$ and $\vartheta \in \tTheta$, and suppose that $\varphi$ and $\vartheta$ are the $F$-orbit and $\tF$-orbit, respectively, of $\xi \in \hat{\bbK}^{\times}$.  Define $\overline{\varphi}$ and $\overline{\vartheta}$ to be the $F$-orbit and $\tF$-orbit, respectively, of $\xi^{-1} = \overline{\xi}$.  Note that $|\varphi| = |\overline{\varphi}|$ and $|\vartheta| = |\overline{\vartheta}|$.  For $\blambda \in \cP^{\Theta}$ and $\bnu \in \cP^{\tTheta}$, define $\overline{\blambda}$ and $\overline{\bnu}$, respectively, by 
 $$ \overline{\blambda}(\varphi) = \blambda(\overline{\varphi}) \; \text{ and } \; \overline{\bnu}(\vartheta) = \bnu(\overline{\vartheta}).$$
 For an element $v \in \Lambda$ or $\tLambda$, we define $\overline{v}$ to be the element of $\Lambda$ or $\tLambda$ obtained when conjugating the coefficients of $v$ when expanding in terms of the $P_{\bmu}$ or $\tP_{\bgamma}$, respectively.
 
 Part (i) of the next result is stated in \cite[1.1.1]{bantan}, where a proof is not given, but is indicated to come from the machinery in \cite[Chapter IV]{Mac} as we have developed it here.  We give the proof of only part (ii) below.
 \begin{lemma} \label{conjchar}
 (i) Let $\blambda \in \cP^{\Theta}$.  Then $\overline{\chi^{\blambda}} = \chi^{\overline{\blambda}}$.\\
 \noindent
 (ii) Let $\bnu \in \cP^{\tTheta}$.  Then $\overline{\psi^{\bnu}} = \psi^{\overline{\bnu}}$.
 \end{lemma}
 \begin{proof} (i): The same as the proof of (ii) below, with appropriate changes made.
 
 \noindent
 (ii):  From Equation (\ref{expand}), it is enough to show that $\overline{\ts_{\bnu}} = \ts_{\overline{\bnu}}$.  Through a series of change of bases, we keep track of what happens to coefficients when expanding $\ts_{\overline{\bnu}}$ in terms of the $\tP_{\bgamma}$.  We have
 $$ \ts_{\overline{\bnu}} = \prod_{\vartheta \in \tTheta} s_{\overline{\bnu}(\vartheta)} (\tY^{(\vartheta)}) = \prod_{\vartheta \in \tTheta} s_{\bnu(\overline{\vartheta})}(\tY^{(\vartheta)}) = \prod_{\vartheta \in \tTheta} s_{\bnu(\vartheta)}(\tY^{(\overline{\vartheta})}).$$
 For a partition $\rho = (\rho_1, \rho_2, \ldots, \rho_{\ell}) \in \cP$, we define the power symmetric function $p_{\rho}$ as $p_{\rho} = p_{\rho_1} p_{\rho_2} \ldots p_{\rho_{\ell}}$.  Recall that the irreducible characters and conjugacy classes of the symmetric group $S_n$ on $n$ letters are both parameterized by partitions of $n$.  Denote the irreducible character of $S_n$ corresponding to the partition $\nu$, where $|\nu| = n$, by $\omega^{\nu}$, and denote the value of this character on an element of cycle-type $\rho$, where $|\rho| = n$, by $\omega^{\nu}(\rho)$.  Recall that all of the values of $\omega^{\nu}$ are integers.  Let $z_{\rho}$ be the size of the centralizer of an element of cycle-type $\rho$.  From the change of basis from Schur functions to power symmetric functions, given in \cite[Proof of I.7.6]{Mac}, we have
 \begin{equation} \label{basis1}
 \ts_{\overline{\bnu}} = \prod_{\vartheta \in \tTheta} \sum_{\rho} \frac{\omega^{\bnu(\vartheta)}(\rho)}{z_{\rho}} p_{\rho_1}(\tY^{(\overline{\vartheta})}) p_{\rho_2}(\tY^{(\overline{\vartheta})}) \cdots p_{\rho_{\ell}}(\tY^{(\overline{\vartheta})}).
 \end{equation}
 Now we change each power symmetric function in the $\tY$ variables to power symmetric functions in the $\tX$ variables using the transform in Equation (\ref{utransform}).  If $\xi \in \vartheta$, then we take $\xi^{-1} = \overline{\xi}$ as the representative of $\overline{\vartheta}$ in the transform, and we make use of the fact that $|\vartheta| = |\overline{\vartheta}|$.  We thus have the change of basis 
 \begin{equation} \label{basis2}
 p_{\rho_i}(\tY^{(\overline{\vartheta})}) = (-1)^{\rho_i|\vartheta| - 1} \sum_{\alpha \in \bar{G}_1^{\tF^{\rho_i |\vartheta|}}} \overline{\xi(\alpha)} p_{\rho_i |\vartheta|/|h_{\alpha}|}(\tX^{(h_{\alpha})}),
 \end{equation}
 where $\alpha \in h_{\alpha}$ and $h_{\alpha} \in \tPhi$.  Finally, we may expand power symmetric functions in the $\tX$ variables in terms of Hall-Littlewood symmetric functions in the $\tX$ variables.  For this change of basis, we have coefficients involving Green's polynomials $Q_{\rho}^{\gamma}(q)$.  It follows from the basis change given in \cite[III.7.1 and 7.8]{Mac} that we have
 \begin{equation} \label{basis3}
 p_k(\tX^{(h)}) = \sum_{\bgamma(h)} Q_{(k)}^{\bgamma(h)}( (-q)^{|h|}) (-q)^{-|h| n(\bgamma(h))}   P_{\bgamma(h)} ( \tX^{(h)}; (-q)^{-|h|}), 
 \end{equation}
 where the coefficients are all rational, by the comment at the beginning of \cite[III.7]{Mac}.  Note that in the change of bases in (\ref{basis1}), (\ref{basis2}), and (\ref{basis3}), the only coefficients when expanding in terms of the $\tP_{\bgamma}$ which are not real occur in (\ref{basis2}), with coefficients of the form $\overline{\xi(\alpha)}$.  If we did the same change of bases to expand $\ts_{\bnu}$, everything would be exactly the same, except these coefficients in (\ref{basis2}) would change to $\xi(\alpha)$.  So, when expanding $\ts_{\overline{\bnu}}$ in terms of $\tP_{\bgamma}$, we obtain the expansion for $\ts_{\bnu}$, except with the coefficients conjugated.  Therefore $\overline{\ts_{\bnu}} = \ts_{\overline{\bnu}}$.
 \end{proof}
 
 It follows immediately from Lemma \ref{conjchar} that an irreducible character $\chi^{\blambda}$ of $G_n$ (or $\psi^{\bnu}$ of $U_n$) is real-valued if and only if $\overline{\blambda} = \blambda$ (or $\overline{\bnu}= \bnu$).  The next result gives a natural bijective correspondence between real-valued irreducible characters of $G_n$ and real-valued irreducible characters of $U_n$ using this combinatorial information.  For $\xi \in \hat{\bbK}^{\times}$, let $[\xi]_{F}$ denote the $F$-orbit of $\xi$ and let $[\xi]_{\tF}$ denote the $\tF$-orbit of $\xi$.
 
\begin{theorem} \label{realcorr}
Let $\blambda \in \cP^{\Theta}_n$ such that $\overline{\blambda} = \blambda$.  Define $r(\blambda) \in \cP^{\tTheta}_n$ by $r(\blambda)([\xi]_{\tF}) = \blambda([\xi]_{F})$.  Then the map $r$ is well-defined, and $\overline{r(\blambda)} = r(\blambda)$.  The map defined by
$$ R: \chi^{\blambda} \mapsto \psi^{r(\blambda)} $$
is a bijection between real-valued irreducible characters of $G_n$ and real-valued irreducible characters of $U_n$.
\end{theorem}
\begin{proof}  We have that $\blambda([\xi]_F) = \blambda([\xi^{-1}]_F)$ for every $\xi \in \hat{\bbK}^{\times}$.  Also, for any $\xi \in \hat{\bbK}^{\times}$ we have
\begin{equation} \label{orbs}
[\xi]_F \cup [\xi^{-1}]_F = [\xi]_{\tF} \cup [\xi^{-1}]_{\tF}.
\end{equation}
We have $r(\blambda)([\xi]_{\tF}) = \blambda([\xi]_F)$ and $r(\blambda)([\xi^{-1}]_{\tF}) = \blambda([\xi^{-1}]_F)$, and it follows from (\ref{orbs}) that $r$ is well-defined and $\overline{r(\blambda)} = r(\blambda)$.

From Lemma \ref{conjchar}, it follows that $R$ maps real-valued irreducible characters of $G_n$ to real-valued irreducible characters of $U_n$.  From (\ref{orbs}), $\blambda$ and $r(\blambda)$ may be viewed as the same partition-valued functions on the unions of orbits $[\xi]_F \cup [\xi^{-1}]_F$, and it follows that $R$ is a bijection.
\end{proof}

\section{Regular characters, semisimple characters, and duality} \label{charduality}

\noindent Let $G$ be a finite group and $N$ a normal subgroup of $G$.  If $\xi$ is a generalized character of $G$, define $T_{G/N}(\xi)$ by
$$ T_{G/N}(\xi) = \frac{1}{|N|} \sum_{n \in N} \xi(ng).$$

Now let $\bar{G}$ be a connected reductive group over $\bar{\bbF}_q$ which is defined over $\bbF_q$ and which has connected center, and let $F$ be a Frobenius map. (We note that this notation slightly conflicts with that used in Section 2 but are confident that no confusion should arise.) Let $W$ be the Weyl group of $\bar{G}$, where $W = \langle s_i \mid i \in I \rangle$, and let $\rho$ be the permutation of the indexing set $I$ which is induced by the action of the Frobenius map $F$.  For any $\rho$-stable subset $J \subseteq I$, let $\bar{P}_J$ be the parabolic subgroup of $\bar{G}$ corresponding to $W_J = \langle s_j \mid j \in J \rangle$, and let $\bar{U}_J$ be the unipotent subgroup.  Let $P_J = \bar{P}_J^F$ and $U_J = \bar{U}_J^F$ be the corresponding parabolic and unipotent subgroups of the finite group $G = \bar{G}^F$.  Define the following operator $*$ on the set of generalized characters of $G$:
$$\xi^* = \sum_{ J \subseteq I \atop {\rho(J) = J}} (-1)^{|J/\rho|} (T_{P_J/U_J}(\xi))^G.$$

As stated in \cite[Chapter 8]{Ca85}, the definition of the operator $*$ and its properties are due to Curtis \cite{curt}, Kawanaka \cite{Ka81}, and Alvis \cite{alvis1, alvis2, alvis3}.  A proof of the following theorem is given in \cite[Section 8.2]{Ca85}.

\begin{theorem} [Curtis, Alvis, Kawanaka]  
The map $\xi \mapsto \xi^*$ is an order 2 isometry of the generalized characters of $G$, so that $\xi^{**} = \xi$ and $\langle \xi, \eta \rangle = \langle \xi^*, \eta^* \rangle$ for all generalized characters $\xi, \eta$ of $G$.
\end{theorem}

The following is immediate from the definition of the map $*$.

\begin{lemma} \label{conj} The map $*$ commutes with complex conjugation.  That is, for any generalized character $\xi$ of $G$, we have $(\overline{\xi})^* = \overline{\xi^*}$.\end{lemma} 

Let $p = {\rm char}(\bbF_q)$, and suppose now that $p$ is a good prime for $\bar{G}$ (see, for example, \cite[Section 1.14]{Ca85} for a definition).  Then we may define a {\em semisimple} character of $G$ to be an irreducible character $\chi$ of $G$ such that $\chi(1)$ is not divisible by $p$.  Recall that the Gelfand-Graev character of $G$, which we will denote by $\Gamma$, is the character of the representation obtained by inducing a non-degenerate linear character from the unipotent subgroup of $G$ up to $G$ (see \cite[Section 8.1]{Ca85} for a full discussion).  A {\em regular} character of $G$ is defined as an irreducible character of $G$ which appears as a constituent of $\Gamma$.  It is well known that the Gelfand-Graev character has a multiplicity free decomposition into irreducible characters of $G$.

The map $*$ gives a duality between the regular characters and semisimple characters of $G$.  The following is proven in \cite[Section 8.3]{Ca85}.

\begin{theorem} \label{duality}
If $\chi$ is a regular character of $G$, then $\chi^* = \pm \psi$, where $\psi$ is a semisimple character of $G$.  If $\psi$ is a semisimple character of $G$, then $\psi^* = \pm \chi$, where $\chi$ is a regular character of $G$.
\end{theorem}

It follows immediately from Lemma \ref{conj} that the duality given by $*$ in Theorem \ref{duality} behaves well when restricted to real-valued characters, as given in the next result.

\begin{corollary} \label{realdual}
The number of real-valued semisimple characters of $G$ is equal to the number of real-valued regular characters of $G$.
\end{corollary}

Now let us restrict our attention to the cases of the finite general linear group, $G_n = {\rm GL}(n, \bbF_q)$, and the finite unitary group, $U_n = {\rm U}(n, \bbF_{q^2})$.  The exact decompositions of the Gelfand-Graev characters of these groups are known, and we give them in terms of the parameterization of characters given in Section \ref{characters}.  For any partition $\lambda \in \cP$, define the {\em length} of $\lambda$, $\ell(\lambda)$, to be the number of non-zero parts of $\lambda$.  Let $\blambda \in \cP^{\cX}$, where $\cX = \Phi, \Theta, \tPhi$, or $\tTheta$.  Define the {\em height} of $\blambda$, written ${\rm ht}(\blambda)$, as 
$$ {\rm ht}(\blambda) = {\rm max}\{ \ell(\blambda(x)) \; \mid \; x \in \cX \}.$$
The decompositions given in the next theorem essentially follow from the more general work of Deligne and Lusztig in \cite{dellusz}, and specific proofs are given for $G_n$ in \cite{zel} and for $U_n$ in \cite{ohm}.
 
\begin{theorem} \label{ggdecomp}
Let $\Gamma$ and $\tilde{\Gamma}$ be the Gelfand-Graev characters of $G_n$ and $U_n$, respectively.  Then the decompositions into irreducibles of $\Gamma$ and $\tilde{\Gamma}$ are
$$ \Gamma = \sum_{\blambda \in \cP^{\Theta}_n \atop { {\rm ht}(\blambda) = 1 }} \chi^{\blambda} \;\;\;\; \text{ and } \;\;\;\; \tilde{\Gamma} = \sum_{\bnu \in \cP^{\tTheta}_n \atop { {\rm ht}(\bnu) = 1 }} \psi^{\bnu}.$$
\end{theorem}

We now count the number of real-valued regular and semisimple characters of the groups of interest.

\begin{theorem} \label{charcount}
Let $G_n = \GL(n, \bbF_q)$ and $U_n = \U(n, \bbF_{q^2})$.  Then:
\begin{align*}
 & \text{the number of real-valued regular characters of $G_n$} \\
=  & \text{the number of real-valued regular characters of $U_n$} \\
= & \text{the number of real-valued semisimple characters of $G_n$} \\
=  & \text{the number of real-valued semisimple characters of $U_n$} \\
= & \left\{\begin{array}{ll} 2q^{m} & \text{ if $q$ is odd and $n = 2m +1$ is odd,} \\ q^{m} + q^{m-1} & \text{ if $q$ is odd and $n = 2m$ is even,} \\ q^{\lfloor n/2 \rfloor} & \text{ if $q$ is even.}\end{array}\right.
\end{align*}
\end{theorem}
\begin{proof} First, from Corollary \ref{realdual}, the number of real-valued regular characters is equal to the number of real-valued semisimple characters in $G_n$ and in $U_n$.  From Lemma \ref{conjchar} and Theorem \ref{ggdecomp}, a real-valued regular character of $G_n$ is of the form $\chi^{\blambda}$, where ${\rm ht}(\blambda) = 1$ and $\overline{\blambda} = \blambda$.  Applying the bijection $R$ given in Theorem \ref{realcorr} to $\chi^{\blambda}$, we obtain some $\psi^{\bnu}$, where $r(\blambda) = \bnu \in \cP^{\tTheta}_n$ satisfies $\overline{\bnu} = \bnu$.  Since the bijection $r$ defined in the proof of Theorem \ref{realcorr} does not change the length of any partition, then $\bnu$ also has the property that ${\rm ht}(\bnu) = 1$.  So, $\psi^{\bnu}$ is a real-valued regular character of $U_n$, and the map $R$ gives a bijection between the real-valued regular characters of $G_n$ and those of $U_n$, and the four quantities of interest are all equal.

It is therefore enough to count the number of $\blambda \in \cP^{\Theta}_n$ such that ${\rm ht}(\blambda) = 1$ and $\overline{\blambda} = \blambda$.  For such a $\blambda$, we have for each $\varphi \in \Theta$ such that $\blambda(\varphi)$ is non-empty, $\blambda(\varphi)$ consists of exactly one part, and $\blambda(\varphi) = \blambda(\overline{\varphi})$.  From the Remark at the beginning of Section \ref{characters}, there is a (non-canonical) bijection between $\Theta$ and $\Phi$ which preserves sizes of orbits, so we may count the $\bmu \in \cP^{\Phi}_n$ such that ${\rm ht}(\bmu) = 1$ and $\overline{\bmu} = \bmu$.  The set $\Phi$ is in bijection with monic irreducible polynomials in $\bbF_q[t]$ with non-zero constant, and so for $f \in \Phi$, we may view $\overline{f}$ as the polynomial in $\bbF_q[t]$ whose roots in $\bar{\bbF}_q$ are the reciprocals of those in $f$.

Now, if $\bmu \in \cP^{\Phi}_n$ and ${\rm ht}(\bmu) = 1$, we may think of $\bmu$ as a collection of monic irreducible polynomials with non-zero constant, $\{ f_i \}$, with a single positive integer $e_i$ associated with each, such that 
$$ \sum_i e_i {\rm deg}(f_i) = n. $$
Furthermore, since $\overline{\bmu} = \bmu$, each $f_i$ satisfies either $\overline{f}_i = f_i$ or the number $e_i$ associated with $f_i$ is equal to that associated with $\overline{f}_i$.  This means that the polynomial
$$f = \prod_i f_i^{e_i}$$ 
satisfies $\overline{f} = f$, where $\overline{f}$ denotes the polynomial in $\bbF_q[t]$ whose roots are the reciprocals of those of $f$.  In fact, such polynomials can always be factored into irreducibles in such a way that each factor $f_i$ either satisfies $\overline{f}_i = f_i$, or $f_i$ occurs with the same power as that of $\overline{f}_i$.  So, the $\bmu \in \cP^{\Phi}_n$ such that ${\rm ht}(\bmu) = 1$ and $\overline{\bmu} = \bmu$ are in one-to-one correspondence with monic polynomials $f \in \bbF_q[t]$ with non-zero constant such that $\overline{f} = f$ and ${\rm deg}(f) = n$.

Let $f \in \bbF_q[t]$, where
$$ f = t^n + a_{n-1} t^{n-1} + \cdots + a_1 t + a_0, \;\;\;\; a_0 \neq 0.$$
Then we have
$$ \overline{f} = a_0^{-1}t^nf(t^{-1}) = t^n + a_0^{-1} a_1 t^{n-1} + \cdots + a_0^{-1} a_{n-1} t + a_0^{-1}.$$
If $f = \overline{f}$, then we must have $a_0^2 = 1$ and $a_{n-i} = a_0^{-1} a_i$ for $1 \leq i < n$.  If $q$ is odd and $n = 2m+1$ is even, then we may choose $a_0 = \pm 1$, and we may choose $a_i$, for $1 \leq i \leq m$, to be any of $q$ elements in $\bbF_q$, and then $a_{n-i}$ must be $ a_0^{-1} a_i$, giving a total of $2q^m$ polynomials.  If $q$ is odd and $n = 2m$ is even, then for $a_0 = 1$, we may choose $a_i$, for $1 \leq i \leq m$, to be any element in $\bbF_q$, while $a_{n-i}= a_i$ for $1 \leq i \leq m-1$, giving $q^m$ polynomials.  If we let $a_0 = -1$, then $a_i$, for $1 \leq i \leq m-1$, may be any element in $\bbF_q$, while $a_m = -a_m$ implies $a_m$ must be $0$, and $a_{n-i} = a_i$ for $1 \leq i \leq m-1$, giving $q^{m-1}$ polynomials, and a total of $q^m + q^{m-1}$.  Finally, if $q$ is even, then we must have $a_0 = 1$, and we may choose $a_i$, for $1 \leq i \leq \lfloor n/2 \rfloor$, to be any of $q$ elements in $\bbF_q$, while the other coefficients are then set, giving $q^{\lfloor n/2 \rfloor}$ polynomials.
\end{proof}

We note that there are formulas for the degrees of irreducible characters of $G_n$ and $U_n$,  and these could be used to count the number of real-valued semisimple characters of these groups directly.  However, giving the decomposition of Gelfand-Graev characters seems to be more straightforward, and the duality given in Corollary \ref{realdual} is very relevant to theme of the main results.

\section{Reality properties and centralizers} \label{realpropsec}

\noindent We let $G_n^+ = G_n \langle \tau \rangle$ and $U_n^+ = U_n \langle \tau \rangle$. All the elements
of the group $G_n^+$ are real, as we noted in \cite{gow1}. The proof of this fact depends critically on a property of the group $G_n$, namely, that all real elements are strongly real, where a strongly real element is one that is inverted by an involution. As we shall see, not all elements of $U_n^+$ are real,
and this non-reality phenomenon is related to the fact that not all real elements of $U_n$ are strongly real. We propose therefore to investigate the question of  whether or not  a real element of $U_n$ is  strongly real. It is the set of unipotent elements that plays the key role in what follows.

\begin{proposition} \label{notstronglyreal} Let $x$ be a regular unipotent element in $U_n$. Then $x$ is not strongly real
if $n$ is even and $q$ is odd, or if $n$ is odd and $q$ is even (note that $x$ is real in all cases).
\end{proposition}

\begin{proof} Let $V$ be a vector space of dimension $n$ over $\bbF_{q^2}$ on which $x$ acts, and let
$f: V\times V\to \bbF_{q^2}$ be a non-degenerate hermitian form preserved by $x$. As we shall prove the results by induction on $n$, we first establish the starting cases. We therefore assume first that $n=2$ and
$q$ is odd. In this case $V$ has a basis consisting of vectors $u$ and $v$ which satisfy
\[
xu=u+v,\quad xv=v.
\]
Since $x$ is an isometry of $f$, we have 
\[
f(u,v)=f(xu,xv)=f(u+v,v)=f(u,v)+f(v,v),
\]
which implies that $f(v,v)=0$. 

Let $s$ be any involution acting on $V$ which inverts $x$. Since $v$ spans the unique one-dimensional
space fixed by $x$, it follows that $sv=\pm v$, and replacing $s$ by $-s$ if necessary, we may assume that $sv=v$. It follows in a straightforward manner that, since $q$ is odd, $su=-u$. Suppose now
that $s$ is an isometry of $f$. Then we have
\[
f(u,v)=f(su,sv)=f(-u,v)=-f(u,v).
\]
Again, since the underlying field has odd characteristic, we obtain $f(u,v)=0$. But the two equalities
$f(u,v)=f(v,v)=0$ imply that $v$ is in the radical of $f$, contradicting the non-degeneracy of $f$.
It follows that $x$ is not strongly real in $U_2$. 

Next, we examine the case that $n=3$ and $q$ is even. $V$ has a basis consisting of vectors
$u$, $v$ and $w$ which satisfy
\[
xu=u+v,\quad xv=v+w,\quad xw=w.
\]
Let $s$ be any involution acting on $V$ which inverts $x$. An elementary calculation, whose details
we omit, shows that 
\[
su=u+av+bw,\quad sv=v+cw,\quad sw=w,
\]
where $a$, $b$ and $c$ are elements of $\bbF_{q^2}$ with either $a=0$, $c=1$ or $a=1$, $c=0$. Now since $x$ is an isometry
of  $f$, the equality $f(v,w)=f(xv,xw)$ yields that $f(w,w)=0$. Similarly, we have
\[
f(u,w)=f(xu,xw)=f(u+v,w)=f(u,w)+f(v,w)
\]
and hence $f(v,w)=0$. We now observe that $f(u,w)\ne 0$, for otherwise we have
\[
f(u,w)=f(v,w)=f(w,w)=0,
\]
which is impossible, since it implies that $w$ is in the radical of $f$. Next, we observe that
\[
f(u,v)=f(xu,xv)=f(u+v,v+w)=f(u,v)+f(u,w)+f(v,v)
\]
and deduce that $f(v,v)=f(u,w)\ne 0$. 

Suppose now that $s$ is an isometry of $f$. Then we must have
\[
f(u,v)=f(su,sv)=f(u+av+bw,v+cw)=f(u,v)+c^qf(u,w)+af(v,v),
\]
since $f(v,w)=f(w,v)=f(w,w)=0$. This implies that $c^qf(u,w)=af(v,v)$. But as we already know
that $ac=0$ and $f(v,v)=f(u,w)\ne 0$, we deduce that $a=c=0$. This contradicts our earlier
observation that one of $a$ and $c$ is 1. Hence, $s$ is not an isometry of $f$, and 
consequently $x$ is not strongly real in $U_3$.

We proceed to the general case where either $n\ge 4$ is even and $q$ is odd or $n\ge 5$ is odd and $q$ is even, and assume that our desired result holds for spaces
of dimension $n-2$. In this case, we can find elements $v$ and $w$ in $V$ with
\[
xv=v+w,\quad xw=w.
\]
A previous argument implies that $f(w,w)=0$. Let $W$ be the one-dimensional subspace of $V$ spanned by $w$ and let $W^{\perp}$ be the subspace of $V$ orthogonal to $W$ (with respect to $f$).
Then $W$ is contained in $W^\perp$, since $f(w,w)=0$. Furthermore, it is a general fact that
$f$ induces a non-degenerate hermitian form $f_1$, say on $W^\perp/W$, which is a space
of dimension $n-2$ over $\bbF_{q^2}$. Since $x$ maps both $W$ and $W^\perp$ into themselves,
it has an induced action as a regular unipotent  element $x_1$, say, on $W^\perp/W$, where it preserves
$f_1$. 

Finally, suppose it is possible that $x$ is inverted by an involutory isometry $s$ of $f$. Then, since
$W$ is the unique one-dimensional subspace of $V$ fixed by $x$, $s$ must also fix $W$ and hence
also leaves $W^\perp$ invariant. We therefore have an induced action of $s$ 
on $W^\perp/W$ as an involutory
isometry, $s_1$ say, of $f_1$, and $s_1$ inverts the regular unipotent element $x_1$. The induction
hypothesis eliminates this possibility and we have a contradiction. It follows that $s$ is not an
isometry of $f$
and $x$ is not strongly real in $U_n$.
\end{proof}

Having shown that certain unipotent elements are not strongly real, we turn to showing that in some sense most real elements of $U_n$ are strongly real. The approach we take is somewhat indirect
and relies on the property of orthogonal groups that all their elements are strongly real. 

We first
recall that there is a one-to-one correspondence between the real classes in $G_n$ and the real
classes in $U_n$ \cite[Theorem 3.8]{gow2}. The correspondence is defined in the following way. Given a real conjugacy  class in $G_n$,
the conjugacy class in $\GL(n,\bbF_{q^2})$ which contains this class is of course real, and, using the Lang-Steinberg theorem, we can show that this conjugacy class in
$\GL(n,\bbF_{q^2})$ intersects $U_n$ in a unique conjugacy class, which is also real. 

Conversely,
given a real conjugacy class in $U_n$, the conjugacy class in $\GL(n,\bbF_{q^2})$ which contains this class is also real, and it intersects $G_n$ in a unique conjugacy class, which is real. Since the real
conjugacy classes of $G_n$ are determined by properties of the elementary divisors of elements,
we can specify unique real conjugacy classes of $U_n$ by the same sets of elementary divisors,
and vice versa.

We also need to make the following observation. Suppose that we have a non-degenerate symmetric
bilinear form of dimension $n$ over $\bbF_q$, where $q$ is odd. Then we may extend this form
to a non-degenerate hermitian form over $\bbF_{q^2}$, and any isometry of the symmetric form
may be extended to an isometry of the hermitian form. Consequently, we can embed any orthogonal
group $\Or(n,\bbF_q)$ into $U_n$. 

Likewise, a non-degenerate alternating bilinear form of dimension
$2m$ over $\bbF_q$, where $q$ is even, may extended to a non-degenerate hermitian form
of the same dimension over $\bbF_{q^2}$, and we may then embed $\hbox{Sp}(2m,\bbF_q)$ into
$U_{2m}$ (this is also true if $q$ is odd, but we will not make use of this embedding).

\begin{proposition} \label{realun} Let $x$ be a real element in $U_n$. Then $x$ is strongly real in either
of the following cases:

\smallskip
\noindent (a) $q$ is odd and each elementary divisor of $x$ of the form $(t\pm 1)^{2m}$ occurs with even multiplicity;

\smallskip
\noindent (b) $q$ is even, $n$ is even, and each elementary divisor of $x$ of the form $(t+ 1)^{2m+1}$ occurs with even multiplicity.
\end{proposition}

\begin{proof} By the earlier discussion, $x$ determines a unique real conjugacy class of $G_n$ with the same elementary divisors. Let $z$ be an element of this conjugacy class of $G_n$.
Then, in case (a), $z$ is conjugate to an element of some orthogonal group $\Or(n,\bbF_q)$
by \cite[p.38, Case C, part (i)]{wall}. Thus we may consider $z$ as an element of 
$\Or(n,\bbF_q)$, and it is  strongly real in this group by a theorem of Wonenburger \cite{won}.
Since we may embed $\Or(n,\bbF_q)$ into $U_n$, we see that $U_n$ has a strongly real
conjugacy class with the same elementary divisors as $z$, and thus the same as those of $x$. This conjugacy class is the same as that of $x$, since the elementary divisors determine the conjugacy class, and hence $x$ is strongly real.

In case (b), we use \cite[p.36, Case B, part (i)]{wall}, together with the fact that 
all elements of $\hbox{Sp}(2m,\bbF_q)$ are strongly real when $q$ is even \cite{ellers} to achieve the desired result.
\end{proof}

We are confident that Proposition \ref{realun} is also a necessary condition for the strong reality of a real
element of $U_n$.  We will not investigate this matter further here, but note that any progress will
involve generalizing considerably the ideas involved in Proposition \ref{notstronglyreal}.

We proceed to examine the reality problem for elements of the coset $U_n\tau$ in the group
$U_n^+$. We begin with some general principles of linear algebra.

Let $x$ be an element of $\GL(n,\bbL)$, where $\bbL$ is an arbitrary field.  We say that $x$ is 
\emph{cyclic} if $x$ acts as a cyclic endomorphism on the underlying $n$-dimensional vector space
over $\bbL$.
We note that if $x$  acts indecomposably on the underlying vector space, it is cyclic, and furthermore,
that $x$ is cyclic 
if and only if its minimal polynomial equals its characteristic polynomial. Finally, if $x$ is cyclic, its centralizer consists of polynomials in $x$. 

It is a theorem of Frobenius \cite[Theorem 66]{kap} that $x$ is conjugate to  $x'$  by a symmetric element, that is, there is an element $s$ of $\GL(n,\bbL)$ with
$s=s'$ and
\[
s^{-1}xs=x'.
\]

We will require the following consequence of Frobenius's theorem.

\begin{lemma} \label{symmetric} Let $x$ be a cyclic element of $\GL(n,\bbL)$ and let $w$ be an element of 
$\GL(n,\bbL)$ which satisfies
\[
w^{-1}xw=x'.
\]
Then $w$ is symmetric.
\end{lemma}

\begin{proof} We know from Frobenius's theorem that a symmetric element $s$ exists satisfying
$s^{-1}xs=x'$. It follows that $w=cs$, where $c$ centralizes $x$. Now 
$c$ is a polynomial in $x$, as $x$ is cyclic, and hence, since $s$ satisfies $s^{-1}xs=x'$, we have
$s^{-1}cs=c'$ also. Finally, we see that
\[
w'=sc'=cs=w,
\]
and thus $w$ is symmetric.
\end{proof}

Let $y\tau$ be an element of $U_n^+$ and let $g=(y\tau)^2$. Then we have $g=
y(y')^{-1}$ and hence
\[
y^{-1}g y=(y')^{-1}y=(g')^{-1}.
\]
Since $g$ and $g'$ are conjugate in the underlying general linear group,
it follows from \cite[p.34, Case A, part (ii)]{wall} that $g$ is a real element of
$U_n$. In the case that $g$ is strongly real, we can prove that $y\tau$ is also strongly
real under suitable hypotheses, as we show below.

\begin{proposition} \label{cyclicstrong} Let $y\tau$ be an element of $U_n^+$ and let $g=(y\tau)^2$.  Suppose that $g$ is cyclic and inverted by an involution $\pi$ in $U_n$.
Then $\pi$ inverts $y\tau$ and hence $y\tau$ is strongly real in $U_n^+$.
\end{proposition}

\begin{proof} As we observed above, 
\[
y^{-1}g y=(y')^{-1}y=(g')^{-1}.
\]
It follows that
\[
y^{-1}\pi^{-1}g \pi y=g'.
\]
Thus, since $g$ is cyclic by hypothesis, Lemma \ref{symmetric} implies that $\pi y$ is symmetric, which translates into
\[
\pi y=y'\pi'.
\]
We now want to prove that
\[
\pi^{-1}(y\tau)\pi=(y\tau)^{-1}=\tau^{-1}y^{-1}=y'\tau.
\]
This amounts to showing that
\[
\pi^{-1}y\pi'=y'.
\]
But this equality holds since $\pi y=y'\pi'$ and $\pi$ is an involution.
Hence $y\tau$ is inverted by $\pi$.
\end{proof}

\begin{corollary} \label{stronglyreal} Let $y\tau$ be an element of $U_n^+$ and let $g=(y\tau)^2$.  Suppose that $g$ is cyclic. Then $y\tau$ is strongly real under either of the following hypotheses:

\smallskip
\noindent(a) $q$ is odd and $g$ has no elementary divisor of the form $(t+ 1)^{2m}$;

\smallskip
\noindent(b) $q$ is even, $n$ is even and $g$ has no elementary divisor of the form
$(t+1)^{2m+1}$.

\end{corollary} 

\begin{proof} As we noted earlier, $g$ is certainly real. Suppose first that $q$ is odd. In this case \cite[Theorem 2.3.1]{wall} shows that any elementary divisor of $g$ of the form $(t-1)^{2m}$
occurs with even multiplicity. Now since $g$ is cyclic by hypothesis, its elementary
divisors occur with multiplicity 0 or 1, and we deduce that $g$ has no elementary divisors
of the form $(t-1)^{2m}$. Furthermore, since $g$ has no elementary divisors of the form
$(t+1)^{2m}$ by hypothesis, $g$ is strongly real by Proposition \ref{realun}, and hence $y\tau$ is strongly
real by Proposition \ref{cyclicstrong}. If we are in case (b),  $g$ is again strongly real by Proposition \ref{realun}
and correspondingly, $y\tau$ is also strongly real by Proposition \ref{cyclicstrong}.
\end{proof}

While Corollary \ref{stronglyreal} gives some information about (strong) reality of elements of $U_n^+$ in the coset
$U_n\tau$, it turns out to be relatively straightforward to show that all elements of $U_n$ are strongly real in $U_n^+$. We begin by proving an analogue for $U_n$ of the theorem of Frobenius
described earlier. Note that in our model of the unitary group $U_n$, if $x\in U_n$, then $x'\in U_n$.

\begin{lemma} \label{unitarysymmetric} Let $x$ be an element of $U_n$. Then there exists a symmetric element
$s$ in $U_n$ with $s^{-1}xs=x'$. 
\end{lemma}

\begin{proof}  We have observed earlier that $x$ and $x'$ are certainly conjugate
in $U_n$.  Suppose first that $x$ is cyclic. Lemma \ref{symmetric} implies that any element which conjugates
$x$ into $x'$ is symmetric, and this proves the lemma in this case. 

In the general case, let
$V$ be the underlying vector space of dimension $n$ over $\bbF_{q^2}$ on which $x$ acts, and let
$f: V\times V\to \bbF_{q^2}$ be a non-degenerate hermitian form preserved by $x$. The results
of Wall show that $V$ is a direct sum of subspaces $V_i$, say, which are 
orthogonal with respect to $f$ and $x$-invariant.
Moreover,  $x$ either acts indecomposably on the subspace $V_i$ or $V_i$ is a direct sum of two
totally isotropic indecomposable $x$-invariant subspaces and the minimal polynomials
of the actions of $x$ on the two summands are relatively prime. 

Let $n_i$ be the dimension of $V_i$
and let $x_i$ be the element of $U_{n_i}$ induced by the action of $x$ on $V_i$. Then $x_i$
is cyclic and hence conjugate to $x_i'$ in $U_{n_i}$ by a symmetric element in this group, by the argument above. It is then straightforward to see that, since $x$ is conjugate in $U_n$ to an
orthogonal direct sum of the $x_i$, it is also conjugate to its transpose by a symmetric element
in $U_n$.
\end{proof}

\begin{corollary} Each element $x$ of $U_n$ is strongly real in $U_n^+$.
\end{corollary}

\begin{proof} Let $s$ be a symmetric element in $U_n$ satisfying $s^{-1}xs=x'$, whose existence is assured by Lemma \ref{unitarysymmetric}. Then we may easily check that $s\tau$ is an involution which inverts $x$.
\end{proof}

We turn to the investigation of some specific elements in $G_n^+$ and $U_n^+$.

\begin{lemma} \label{unipotentsquare} Suppose that $q$ is odd. Then the following hold.

\smallskip
\noindent (a) 
If $n$ is odd, then there is an element $x \in G_n$ such that $(x\tau)^2 $ is regular unipotent, and there is an element $y \in U_n$ such that $(y\tau)^2$ is regular unipotent.  The element $y\tau$ is 
strongly real in $U_n^+$.  

\smallskip
\noindent(b)  If $n$ is even, there there is an element $x \in G_n$ such that $(x\tau)^2 = -u$, where $u \in G_n$ is regular unipotent, and there is an element $y\in U_n$ such that $(y\tau)^2 = -v$, where $v \in U_n$ is regular unipotent.
\end{lemma}

\begin{proof} We consider case (a) first. Since a regular unipotent element in $G_n$ has the single
elementary divisor $(t-1)^n$, it follows from \cite[Theorem 2.3.1]{wall} that there is an element $x \in G_n$ such that $(x\tau)^2 $ is regular unipotent (note that, in Wall's terminology, multipliers in this context are precisely elements of the form $(x\tau)^2$). Now let $\phi$ be the map of conjugacy classes described in Lemma \ref{correspondence} and let $\phi[x\tau]=[y\tau]$. Lemma \ref{conjugacy} shows that $(x\tau)^{-2}$ and
$(y\tau)^2$ are conjugate in $\GL(n,\bbK)$, and this implies that $(y\tau)^2$ is also regular unipotent.
That $y\tau$ is strongly real follows from Corollary \ref{stronglyreal}.

The proof in case (b) is similar, since if $u$ is regular unipotent, $-u$ has the single elementary divisor $(t+1)^n$ and  hence is a multiplier by Wall's theorem. This implies the existence
$x \in G_n$ such that $(x\tau)^2 = -u$.  Lemma \ref{conjugacy} then implies the  existence
$y \in U_n$ such that $(y\tau)^2 = -v$, where $v$ is regular unipotent in $U_n$.
\end{proof} 

We note that when $q$ is even and when $n$ is odd, there are also elements $x\tau \in G_n^+$ and $y\tau \in U_n^+$ such that $(x\tau)^2$ and $(y\tau)^2$ are regular unipotent in $G_n$ and $U_n$, respectively.  Such elements are explicitly given in Section \ref{chartwosec}.

\begin{lemma} \label{centorder} Let $n = 2m+1$ and $q$ both be odd, and $x \in G_n$ be an element such that $(x\tau)^2$ is regular unipotent.  Then
\[
 |C_{G_n^+}(x\tau)| = 4q^m.
\]
\end{lemma}

\begin{proof} Clearly, we have
\[
|C_{G_n^+}(x\tau)| = 2|C_{G_n}(x\tau)|.
\]
We set $u=(x\tau)^2$ and note then that $ux'=x$. We observed in the introduction that the centralizer
of $x\tau$ in $G_n$ is identical with the isometry group of the the bilinear form, $b$ say, defined by $x$.
Now in our case $x+x'=(u+1)x'$ is symmetric and invertible, since $u+1$ is invertible, and hence
determines a non-degenerate symmetric bilinear form, $f$ say. Furthermore,
$u$ is an isometry of $f$.  Fulman and Guralnick observe in \cite[p.386]{fulman} that the isometry
group of $b$ is identical with the centralizer of $u$ in the isometry group of $f$, which in this case is the
orthogonal group $\Or(2m+1, \bbF_q)$. Since by \cite[p.38, Case C, part (iv)]{wall}, this centralizer
has order $2q^m$, we obtain the desired result.
\end{proof}

\begin{corollary} \label{centunitary} Let $n=2m+1$ and $q$ both be odd, and let  $y$ be an element in $U_n$ such that $(y\tau)^2$ is regular unipotent. Then we have
$$ |C_{U_n^+}(y\tau)| = 4q^m.$$
\end{corollary}

\begin{proof} We know that $y$ exists from Lemma \ref{unipotentsquare} and the formula for the order of the centralizer follows from  Lemma \ref{centralizers} and Lemma \ref{centorder}.
\end{proof}

\section{Character values} \label{charvaluessec}

\noindent Let $\chi$ is an irreducible complex character of a finite group $G$.  Recall that the Frobenius-Schur indicator, $\varepsilon(\chi)$, takes the value $1$ if $\chi$ may be realized over the real field, $-1$ if $\chi$ is real-valued but  cannot be realized over the real field, and $0$ if $\chi$ is not real-valued.  The following result on Frobenius-Schur indicators is proven in \cite{gow1} for the case of $G_n$ and $G_n^+$, and in \cite{TV} in the cases for $U_n$ and $U_n^+$.

\begin{theorem} \label{indicators} Let $G_n^+$ and $U_n^+$ be the split extensions of $G_n$ and $U_n$ by the transpose inverse automorphism, respectively.  We have the following:
\begin{enumerate} 
\item Let $\theta$ be a real-valued character of $G_n$.  Then $\varepsilon(\theta) = 1$, and $\theta$ has two extensions $\chi$ and $\chi'$ to $G_n^+$ such that $\varepsilon(\chi) = \varepsilon(\chi') = 1$.
\item Let $\theta$ be a character of $U_n$ such that $\varepsilon(\theta) = 1$.  Then $\theta$ has two extensions $\chi$ and $\chi'$ to $U_n^+$ such that $\varepsilon(\chi) = \varepsilon(\chi') = 1$.  
\item Let $\theta$ be a character of $U_n$ such that $\varepsilon(\theta) = -1$.  Then $\theta$ has two extensions $\chi$ and $\chi'$ to $U_n^+$ such that $\varepsilon(\chi) = \varepsilon(\chi') = 0$.
\end{enumerate}
\end{theorem} 

Except possibly when $n$ and $q$ are both even, the group $U_n$ has irreducible characters $\theta$ such that $\varepsilon(\theta) = -1$, and Theorem \ref{indicators}, part (3),  above says that when we extend such $\theta$ to $U_n^+$, some of the character values will not be real.  The following result tells us that these values  are purely imaginary (and all
other values are 0).

\begin{lemma} \label{extlemma} Let $N$ be a normal subgroup of index 2 of some finite group $G$.  Let $\theta$ be a real-valued irreducible character of $N$ which is invariant under conjugation by elements in $G$.  Let $\chi$ be an extension of $\theta$ to $G$, and suppose that $\chi$ is not real-valued.  Then for every $g \in G\setminus N$, $\chi(g)$ is either $0$ or purely imaginary.
\end{lemma}
\begin{proof} We have $\theta^G = \chi + \sigma \chi$, where $\sigma$ is the sign character of $G/N$.  Since $\theta^G$ is real-valued, $\bar{\chi}$ is a constituent of $\theta^G$, but $\chi \neq \bar{\chi}$ since $\chi$ is not real-valued, so $\bar{\chi} = \sigma \chi$.  For $g \in G \setminus N$, $\sigma(g) = -1$, and so $\bar{\chi}(g) = -\chi(g)$, hence $\chi(g) = 0$ or is purely imaginary.  \end{proof}

We may immediately apply Lemma \ref{extlemma} to see that the real-valued irreducible characters of 
Frobenius-Schur indicator $-1$ of $U_n$ vanish on 
many elements of $U_n\tau$ when they are extended to $U_n^+$.

\begin{corollary} \label{purelyimaginary} Let $\theta$ be an irreducible character of $U_n$ such that $\varepsilon(\theta) = -1$, and let $\chi$ be an extension of $\theta$ to $U_n^+$.  Then $\chi(g\tau) = 0$ for any real element
$g\tau$ of the coset $U_n\tau$ in $U_n^+$. 
\end{corollary}

\begin{proof} Certainly, $\theta$ is real-valued, but from Theorem \ref{indicators}, part (3), $\chi$ is not real-valued.  From Lemma \ref{extlemma}, $\chi(g\tau)$ is either $0$ or purely imaginary.  But $g\tau$ is a real element by hypothesis and so $\chi(g\tau)$ must be a real number. This implies that $\chi(g\tau)=0$. 
\end{proof} 

As in Section \ref{charduality}, let $\bar{G}$ be a connected reductive group over $\bar{\bbF}_q$ which is defined over $\bbF_q$ and which has connected center, let $F$ be a Frobenius map, and let $G = \bar{G}^F$.  Also assume that $p = {\rm char}(\bbF_q)$ is a good prime for $\bar{G}$, and recall that a semisimple character of $G$ is an irreducible character with degree prime to $p$.  We note that if $\bar{G} = {\rm GL}(n, \bar{\bbF}_q)$, then every prime is a good prime for $\bar{G}$.  Green, Lehrer, and Lusztig \cite{GLL} found that the character values of $G$ on a regular unipotent element of $G$ can only be $0$, $1$, or $-1$, and are congruent to the character degree modulo $p$, as stated in the next result.

\begin{theorem} [Green, Lehrer, Lusztig] \label{GLL} Let $\chi$ be an irreducible character of $G$, let $u \in G$ be a regular unipotent element, and let $p = {\rm char}(\bbF_q)$.  If $\chi(1)$ is prime to $p$, then $\chi(u) = \pm 1$, and otherwise $\chi(u) = 0$.  Also,
$$\chi(1) \equiv \chi(u) \pmod  {p}.$$
\end{theorem}

Our main result in Theorem \ref{extensionvalues} below may be viewed as a generalization of Theorem \ref{GLL} to the groups $G_n^+$ and $U_n^+$.  Before giving this result we first prove the following.

\begin{lemma} \label{galois}
Let $\chi$ be a character of $G_n^+$ or $U_n^+$ which is an extension of a real-valued semisimple character of $G_n$ or $U_n$, respectively.  Let $\alpha \in {\rm Gal}(\bar{\mathbb{Q}}/\mathbb{Q})$.  Then $\chi^{\alpha} = \alpha \circ \chi$ is a character of $G_n^+$ or $U_n^+$ which is an extension of a real-valued semisimple character of $G_n$ or $U_n$, respectively.
\end{lemma}
\begin{proof} The proof is the same in either case.  Let $\chi$ be such a character of $G_n^+$, and let $\chi|_{G_n}$ denote restriction to $G_n$.  First note that we have
$$(\chi^{\alpha})|_{G_n} = (\chi|_{G_n})^{\alpha}.$$
Now, since $\chi|_{G_n}$ is an irreducible character, so is $(\chi|_{G_n})^{\alpha}$, and so $(\chi^{\alpha})|_{G_n}$ is irreducible.  This implies that $\chi^{\alpha}$ is the extension of a real-valued character.  Since $\chi(1)$ is prime to $p$, then so is $\chi^{\alpha}(1) = \chi(1)$, and the result follows.
\end{proof}

\begin{theorem} \label{extensionvalues} Let $n = 2m+1$ and $q$ both be odd.  Let $y\tau$ be an element of $G_n^+$ (or $U_n^+$) such that $(y\tau)^2$ is regular unipotent, and let $\chi$ be a character of $G_n^+$ (or $U_n^+$) which is an extension of a real-valued irreducible character of $G_n$ (or $U_n$).  Then $\chi(y\tau) = \pm 1$ if $\chi(1)$ is prime to $p$, and $\chi(y\tau) = 0$ otherwise.  Also, for any irreducible $\chi$ of $G_n^+$ (or $U_n^+$),
$$ \chi(\tau) \equiv \pm\chi(y\tau) \pmod {p}.$$
\end{theorem}
\begin{proof} We give the proof in the case of $U_n$ and $U_n^+$. The proof for
$G_n$ and $G_n^+$ is identical.
Let $g = y\tau \in U_n^+$, so $g^2 = u \in U_n$ is regular unipotent.  Let $\chi_1, \chi_2, \ldots, \chi_t$, be the irreducible characters of $U_n^+$ which are extended from real-valued semisimple characters of $U_n$.  From Theorems \ref{charcount} and \ref{indicators}, we have $t = 4q^m$, where $n = 2m+1$.  Let $\mathcal{O}$ be the ring of algebraic integers in $\bar{\mathbb{Q}}$.  Let $\chi$ be one of the $\chi_j$.  We have
$$ \chi(g)^2 \equiv \chi(g^2) \pmod  {2\mathcal{O}}, $$
and from Theorem \ref{GLL}, $\chi(g^2) = \chi(u) = \pm 1$.  So
\begin{equation} \label{mod2}
\chi(g)^2 \equiv 1 \pmod  {2\mathcal{O}}, 
\end{equation}
and in particular, $\chi(g) \neq 0$.  Now, from the arithmetic-geometric mean inequality, we have
\begin{equation} \label{agm}
\frac{1}{t} \sum_{j = 1}^t |\chi_j(g)|^2 \geq \prod_{j=1}^{t} |\chi_j(g)|^{2/t}.\end{equation}
From Lemma \ref{galois}, the set $\{ \chi_1, \ldots \chi_t \}$ is a union of Galois orbits, and since $\chi_j(g)$ is an algebraic integer, we have
$$\prod_{j=1}^t \chi_j(g) \in \mathbb{Z}, \;\;\; \text{ and so } \;\;\; \prod_{j=1}^t |\chi_j(g)|^{2/t} \geq 1.$$
From (\ref{agm}), it follows that
\begin{equation} \label{ineq1}
\sum_{j = 1}^t |\chi_j(g)|^2 \geq t = 4q^m.
\end{equation}
Corollary  \ref{centunitary} tells us that $|C_{U_n^+}(g)| = 4q^m$, and so by the column orthogonality relation of characters, we have
\begin{equation} \label{ineq2}
\sum_{j = 1}^t |\chi_j(g)|^2 \leq 4q^m.
\end{equation}
The inequalities (\ref{ineq1}) and (\ref{ineq2}) together tell us that for every $j$, we have $|\chi_j(g)| = 1$, and if $\chi$ is any other irreducible character of $U_n^+$, then $\chi(g) = 0$.  From Lemma \ref{unipotentsquare}, $\chi_j(g)$ is real, and so we have $\chi_j(g) = \pm 1$.

For the second statement, first consider the case $G_n$.  Since $(y\tau)^2 = u$ is regular unipotent, then $\big((y\tau)^{p^k}\big)^2 = I$ for some $k$, where $I$ is the identity matrix.  This implies that $(y\tau)^{p^k} = s\tau$ for some symmetric matrix $s \in G_n$, since for any $g \in G_n$, $gs\tau g^{-1} = gsg'\tau$.  So, from the classification of symmetric bilinear forms, and since $n$ is odd, we have $s\tau$ is conjugate to either $\tau$ or $dI\tau$, where $d$ is some non-square in $\bbF_q$.  If $(y\tau)^{p^k}$ is conjugate to $\tau$, then we have $\chi(\tau) \equiv \chi(y\tau) \pmod {p}$ for any irreducible $\chi$ of $G_n^+$.  If $(y\tau)^{p^k}$ is conjugate to $dI\tau$, let $\chi$ be an irreducible character of $G_n^+$ which is extended from a real-valued irreducible of $G_n$ (otherwise, $\chi(y\tau) = \chi(\tau) = 0$).  Let $\Pi$ be the representation of $G_n^+$ with character $\chi$.  Then $\Pi(dI\tau) = \Pi(dI)\Pi(\tau)$, and since $\Pi$ restricted to $G_n$ has real-valued character, then its central character on $G_n$ takes only the values $\pm 1$.  Since $dI$ is in the center of $G_n$, then $\Pi(dI)\Pi(\tau) = \pm\Pi(\tau)$.  So, $\chi(dI\tau) = \pm \chi(\tau)$.  Now we have $\chi(\tau) \equiv \pm \chi(y\tau) \pmod {p}$.  

In the case $U_n$, again we have $(y\tau)^{p^k} = s\tau$ for some $k$, and for some symmetric $s$ in $U_n$.  The conjugacy classes in $U_n \tau$ of order $2$ are again in correspondence with $U_n$-equivalence classes of symmetric matrices in $U_n$, and by Theorem \ref{classcorrespondence} there are exactly two such classes, since there are two such classes in $G_n$.  It follows from the classification of symmetric bilinear forms, the fact that $n$ is odd, and the definition of $U_n$, that these two classes are represented by $I$ and $bI$, where $b$ is an element of $ M= \{ a \in \bbF_{q^2} \; | \; a^{q+1} = 1 \}$ which is not the square of an element of $M$.  Thus, $s\tau$ is conjugate to either $\tau$ or $bI\tau$, and by the same argument as above, for any irreducible $\chi$ of $U_n^+$, we have $\chi(\tau) \equiv \pm \chi(y\tau) \pmod {p}$.
\end{proof}

The next result follows directly from Theorems \ref{extensionvalues} and \ref{GLL}.

\begin{corollary} Let $\chi$ be an irreducible character of $G_n^+$ (or $U_n^+$) which is an extension of a real-valued irreducible of $G_n$ (or $U_n$).  Then
$$ \chi(1) - \chi(\tau) \equiv \pm 2 \text{ or } 0 \pmod {p}.$$
\end{corollary} 

We now show that there is no direct analogue of Theorem \ref{extensionvalues} in the even dimensional case for $U_n$.  Suppose that $q$ is odd and $n$ is even, and let $u$ be a regular unipotent
element in $U_n$. We showed in Proposition  \ref{notstronglyreal} that $u$ is real but not strongly real in $U_n$. We also noted in \cite{gow3}, discussion before Theorem 4.4,  that there exists a real-valued irreducible character
$\theta$ of $U_n$ with $\varepsilon(\theta)=-1$ and $\theta(u)\ne 0$. It follows from Theorem \ref{GLL}
that $p$ does not divide $\theta(1)$ and $\theta(u)=\pm 1$. 

By Theorem \ref{indicators}, $\theta$ extends to an irreducible character $\chi$ of $U_n^+$ which is not real-valued. We also know by Lemma \ref{unipotentsquare} that there is an element
$y\in U_n$ with $(y\tau)^2=-u$. We show next that $\chi(y\tau)$ is a non-zero purely imaginary
complex number.

\begin{theorem} \label{complexextensionvalues} Let $n$ be even and $q$ odd, and let $u$ be
a regular unipotent element in $U_n$.   
Let $y\tau$ be an element of  $U_n^+$ such that $(y\tau)^2=-u$ and let  $\theta$ be a 
real-valued irreducible character of  $U_n$ of degree prime to $p$ for which $\varepsilon(\theta)=-1$.
Let $\chi$ 
be an extension of $\theta$ to an irreducible character of  $U_n^+$.  Then $\chi(y\tau) $ is a non-zero
purely imaginary complex number. Hence $y\tau$ is not real.
\end{theorem}

\begin{proof} We first observe that our remarks above show that there exist characters $\theta$
with the stated property.
Let $g = y\tau \in U_n^+$, so that $g^2 = -u \in U_n$ is regular unipotent.   Let $\mathcal{O}$ be the ring of algebraic integers in $\bar{\mathbb{Q}}$. As in the proof of Theorem \ref{extensionvalues},  we have
$$
 \chi(g)^2 \equiv \chi(g^2) \pmod  {2\mathcal{O}},
  $$
and from Theorem \ref{GLL}, 
\[
\chi(g^2) = \chi(-u) = \pm \theta(u)=\pm 1.
\]
Hence
\[
\chi(g)^2 \equiv 1 \pmod  {2\mathcal{O}}, 
\]
and in particular, $\chi(g) \neq 0$.  It follows from Corollary \ref{purelyimaginary} that
$\chi(g)$ is a non-zero purely imaginary complex number and hence $g=y\tau$ is not real.
\end{proof}

On the basis of examining examples, we conjecture that if $n=2m$, there are $q^{m-1}$ characters
$\theta$ satisfying the hypothesis of Theorem \ref{complexextensionvalues}, and if
$\chi$ is an extension of $\theta$, then $\chi(y\tau)=\pm \sqrt{-q}$.  Examples also suggest that there should be corresponding irreducible characters $\psi$ of $G_n^+$, extended from real-valued characters of $G_n$, such that $\psi(x\tau) = \pm \sqrt{q}$, where $(x\tau)^2 = -u$ and $u$ is regular unipotent in $G_n$.

\section{Characteristic Two} \label{chartwosec}

In this section we apply the theory of Gelfand-Graev characters on disconnected reductive groups, due to Sorlin \cite{sorlin1, sorlin2}, to obtain the results on extended character values when our finite field has characteristic 2.  The reference \cite{sorlin1} is a summary of the main results of the theory of Gelfand-Graev characters of disconnected groups, while \cite{sorlin2} contains all proofs for the statements.  We first establish that our particular example fits the general framework of the theory of Sorlin.  

Let $\bar{G}_n = {\rm GL}(n, \bar{\bbF}_q)$, where $q$ is a power of $2$.  Define the automorphism $\sigma$ on $\bar{G}$ by $\sigma(g) = w_0 (g')^{-1} w_0$, where $g'$ denotes transpose as before, and $w_0$ is the element with $1$'s on the antidiagonal and $0$'s elsewhere.  We let $F$ be the standard Frobenius map, but now we define $\tilde{F}$ by $\tilde{F} = F \circ \sigma$.  Note that this Frobenius map differs from the $\tilde{F}$ defined earlier as it includes conjugation by $w_0$, but it follows from the Lang-Steinberg Theorem that $\bar{G}_n^{\tilde{F}}$ is isomorphic to the finite unitary group $U_n$ as defined before.  The automorphism $\sigma$ commutes with $F$ and $\tilde{F}$, and so $\sigma$ is {\em rational} with respect to these Frobenius maps.  We let $\bar{G}_n\langle \sigma \rangle$ denote the semidirect product of $\bar{G}$ by $\sigma$, making $\bar{G}_n\langle \sigma \rangle$ a disconnected reductive group.  Note that if $\tau$ is the transpose inverse automorphism, then since $w_0 \in \bar{G}_n$, we have $\bar{G}_n\langle \sigma \rangle = \bar{G}_n\langle \tau \rangle$.  The unipotent elements of this group are described in \cite[I.2.7]{Spalt}, and in particular, $\sigma$ is unipotent.  Also note that $\sigma$ stabilizes the standard Borel subgroup and its maximal torus in $\bar{G}$, since we have conjugated by $w_0$.  It follows from \cite[Cor. 1.33]{dignemichel} that $\sigma$ is a rational {\em quasi-central} automorphism of $\bar{G}$.

The theory in \cite{sorlin1, sorlin2} now allows us to define Gelfand-Graev characters of the groups $\bar{G}_n^{F} \langle \sigma \rangle$ and $\bar{G}_n^{\tilde{F}} \langle \sigma \rangle$.  If $\bar{N}$ is the standard unipotent subgroup of $\bar{G}_n$, we note that $\bar{N}$ is fixed under $F$, $\tilde{F}$, and $\sigma$.  Let $N$ denote either $\bar{N}^{F}$ or $\bar{N}^{\tilde{F}}$, which are the standard unipotent subgroups of $\bar{G}_n^{F} = G_n$ and $\bar{G}_n^{\tilde{F}} \cong U_n$, respectively, and let $G$ denote either $\bar{G}_n^{F}$ or $\bar{G}_n^{\tilde{F}}$.  Since ${\rm char}(\bbF_q)=2$, there exist $\sigma$-fixed non-degenerate linear characters of $N$, and these are exactly the real-valued non-degenerate linear characters of $N$.  Choose one of these characters $\theta$, and extend it to a linear character $\theta^+$ of $N\langle \sigma \rangle$ such that $\theta^+(\sigma) = 1$, which is possible since $\theta$ is real-valued.  The {\em Gelfand-Graev character} $\Gamma$ of $G\langle \sigma \rangle$ is defined as (see \cite[Prop. 5.1]{sorlin2}) 
$$ \Gamma = {\rm Ind}_{N\langle \sigma \rangle}^{G\langle \sigma \rangle}(\theta^+).$$
Now note that if $n = 2m$ is even, then $\bar{G}^{\sigma} \cong {\rm Sp}(2m, \bar{\bbF}_q)$, and if $n = 2m + 1$ is odd, then $\bar{G}^{\sigma} \cong {\rm O}(2m+1, \bar{\bbF}_q) \cong {\rm Sp}(2m, \bar{\bbF}_q)$.  Then $\bar{G}^{\sigma}$ has semisimple rank $m = \lfloor n/2 \rfloor$, and $Z(\bar{G}^{\sigma})$ consists of only the identity, and in particular is connected.  From \cite[Cor. 8.12]{sorlin2}, it follows that there is a unique Gelfand-Graev character of $G\langle \sigma \rangle$, and so $\Gamma$ does not depend on the linear character $\theta$.  From \cite[Prop. 6.1]{sorlin2}, the character $\Gamma$ of $G\langle \sigma \rangle$ is multiplicity-free, and its restriction to $G$ is exactly the Gelfand-Graev character $\Gamma_G$ of $G$.

Consider the set of complex-valued $G$-class functions on the coset $G\sigma$.  Define an inner product on such functions by 
\begin{equation} \label{inner}
\langle \alpha, \beta \rangle_{G\sigma} = \frac{1}{|G|}\sum_{x\sigma \in G\sigma} \alpha(x\sigma) \overline{\beta(x\sigma)}.
\end{equation}
A duality operation on these class functions, much like the duality discussed in Section \ref{charduality} of this paper, is defined in \cite[Def. 3.10]{dignemichel}, and we will use the notation $\alpha^*$ for the dual of $\alpha$.  By \cite[Cor. 3.12]{dignemichel}, the operation $*$ is an isometric involution with respect to the inner product (\ref{inner}).  Denote by $\Gamma_{G\sigma}$ the Gelfand-Graev character of $G\langle\sigma\rangle$ restricted to $G\sigma$, and let $\Xi_{G\sigma} = \Gamma_{G\sigma}^*$.  Let $\langle \cdot, \cdot \rangle$ and $\langle \cdot, \cdot \rangle_G$ denote the standard inner products on class functions of $G\langle\sigma\rangle$ and $G$, respectively.     

\begin{lemma} \label{duallemma} Let $\chi$ be an irreducible character of $G\langle \sigma \rangle$ extended from a $\sigma$-stable irreducible character of $G$, and let $\chi_{G\sigma}$ be the restriction of $\chi$ to $G\sigma$.  Then we have
$$ \langle \chi_{G\sigma}, \Xi_{G\sigma} \rangle_{G\sigma} = \pm 1 \text{ or } 0.$$
\end{lemma}
\begin{proof}   We have 
$$\langle \chi, \chi \rangle = \langle \chi_G, \chi_G \rangle_G = \langle \chi_{G\sigma}, \chi_{G\sigma} \rangle_{G\sigma} = 1.$$  
Define $\chi^*$ to be the class function on $G\langle \sigma \rangle$ such that $\chi^*$ restricted to $G\sigma$ is $\chi_{G\sigma}^*$, and $\chi^*$ restricted to $G$ is $\chi_G^*$, where the latter $*$ denotes the duality operation defined in Section \ref{charduality}.  We have
$$\langle \chi^*, \chi^* \rangle = \frac{1}{2} \langle \chi_G^*, \chi_G^* \rangle_G + \frac{1}{2} \langle \chi_{G\sigma}^*, \chi_{G\sigma}^* \rangle_{G\sigma} = 1,$$
and so $\pm \chi^*$ is an irreducible character of $G\langle \sigma \rangle$.  We have
$$ \langle \chi^*, \Gamma \rangle = \langle \chi^*_G, \Gamma_G \rangle_G = \pm 1 \text{ or } 0,$$
since both $\Gamma$ and $\Gamma_G$ are multiplicity free.  Finally, we have
$$ \langle \chi_{G\sigma}, \Xi_{G\sigma} \rangle_{G\sigma} = \langle \chi_{G\sigma}^*, \Gamma_{G\sigma} \rangle = 2\langle \chi^*, \Gamma \rangle - \langle \chi^*_G, \Gamma_G \rangle_G = \pm 1 \text{ or } 0.    \qedhere$$
\end{proof}

It follows from \cite[Thm. 8.4(ii)]{sorlin2} and the fact that $*$ is an isometry that
$$\langle \Gamma_{G\sigma}, \Gamma_{G\sigma} \rangle_{G\sigma} =  \langle \Xi_{G\sigma}, \Xi_{G\sigma} \rangle_{G\sigma} = q^{\lfloor n/2 \rfloor}.$$
We note that this also follows from our Theorem \ref{charcount}, as we now explain.  If $\psi$ is an irreducible character of $G$ which is a constituent of $\Gamma_G$, consider $\psi$ induced to $G\langle \sigma \rangle$, which we denote by $\psi^{G\langle \sigma \rangle}$.  We have
$$ \langle \Gamma, \psi^{G\langle \sigma \rangle} \rangle = \langle \Gamma_G, \psi \rangle_G = 1.$$
So, if $\psi$ is a real-valued character, then $\psi^{G\langle\sigma\rangle} = \psi_1+ \psi_2$, where $\psi_1$ and $\psi_2$ are irreducible extensions of $\psi$ to $G\langle \sigma \rangle$, and so exactly one of these extensions of $\psi$ to $G\langle \sigma \rangle$ is a constituent of $\Gamma$.  If $\psi$ is not real-valued, then $\overline{\psi}$ is also a constituent of $\Gamma$ since $\Gamma$ is real-valued, and $\psi^{G\langle \sigma \rangle}$ is an irreducible character of $G\langle \sigma \rangle$ which is $\psi + \overline{\psi}$ on $G$ and $0$ on $G\sigma$.  Thus, $\langle \Gamma, \Gamma \rangle$ is equal to the number of real-valued constituents of $\Gamma_G$ plus half the number of constituents of $\Gamma_G$ which are not real-valued.  From Theorem \ref{charcount}, we now have
$$ \langle \Gamma_{G\sigma}, \Gamma_{G\sigma} \rangle_{G\sigma} = 2\langle \Gamma, \Gamma \rangle - \langle \Gamma_G, \Gamma_G \rangle_G = q^{\lfloor n/2 \rfloor}.$$

The definition of a {\em regular unipotent} element in a disconnected reductive group is given in \cite[I.4.8]{Spalt}, and it follows from that section that all regular unipotent elements of $\bar{G}_n \langle \sigma \rangle$ are conjugate. From \cite[Prop. II.10.2]{Spalt}, all regular unipotent elements of $\bar{G}_n \langle \sigma \rangle$ are of the form $v\sigma$ with $v \in \bar{G}_n$.  A {\em rational} regular unipotent element of $\bar{G}_n^F\langle \sigma \rangle$ (or $\bar{G}_n^{\tilde{F}}\langle \sigma \rangle$) is a regular unipotent element of $\bar{G}_n \langle \sigma \rangle$ which is fixed under the Frobenius map $F$ (or $\tilde{F}$).  The result \cite[Prop. II.10.2]{Spalt} can be used to give explicit regular unipotent elements in $\bar{G}_n\sigma$.  If $m = \lfloor n/2 \rfloor$, then $v\sigma$ is regular unipotent, where $v_{ii} = 1$ for $1 \leq i \leq n$, $v_{i, i+1} = 1$ for $1 \leq i \leq m$, and $v_{ij}=0$ otherwise.  These are in fact rational regular unipotent elements in $\bar{G}_n^{F}\sigma$, where if $n$ is odd, $(v\sigma)^2$ is regular unipotent in $G_n$, and if $n$ is even, $(v\sigma)^2$ is unipotent of type $(n-1,1)$ in $G_n$.  It follows that all of the rational regular unipotent elements in $\bar{G}_n^{F}\sigma$ are conjugate.  Applying Theorem \ref{classcorrespondence}, we obtain exactly one conjugacy class of rational regular unipotent elements in $\bar{G}_n^{\tilde{F}}\sigma$ as well.

It follows from \cite[Thm. 3.6]{sorlin2} that 
\begin{equation} \label{unicount}
\text{the number of rational regular unipotent elements in $G\sigma$ is } \frac{|G|}{q^{\lfloor n/2 \rfloor}}.
\end{equation}
We note that given the description of rational regular unipotent elements above, we can also count the number of such elements in $\bar{G}_n^{F}\sigma$ using \cite[Thm. 3.3]{fulman} and apply Theorem \ref{classcorrespondence} to find that this is also the number of such elements in $\bar{G}_n^{\tilde{F}}\sigma$.

\begin{proposition} \label{regunipotentvalues}
Let $G$ denote either $\bar{G}_n^{F} = G_n$ or $\bar{G}_n^{\tilde{F}} \cong U_n$.  Let $v\sigma$ be a regular unipotent element in $G\sigma$.  For any irreducible character $\chi$ of $G\langle \sigma \rangle$, $\chi(v\sigma) = \pm 1$ or $0$.
\end{proposition}
\begin{proof} From \cite[Cor. 8.12]{sorlin2}, $\Xi_{G\sigma}$ takes the value $q^{\lfloor n/2 \rfloor}$ on regular unipotent elements and $0$ elsewhere.  Let $\chi$ be any irreducible of $G\langle \sigma \rangle$.  Applying (\ref{unicount}) and the fact that the rational regular unipotent elements in $G\sigma$ are all conjugate, we have
\begin{align*}
\langle \chi_{G\sigma}, \Xi_{G\sigma} \rangle_{G\sigma}  &= \frac{1}{|G|} \sum_{v\sigma \text{ regular} \atop{\text{unipotent}}} \chi(v\sigma) \Xi_{G\sigma}(v\sigma)\\
&= \frac{1}{|G|} \frac{|G|}{q^{\lfloor n/2 \rfloor}} \chi(v\sigma) q^{\lfloor n/2 \rfloor} = \chi(v\sigma).
\end{align*}
From Lemma \ref{duallemma}, we have $\chi(v\sigma) = \pm 1$ or $0$.
\end{proof}

Finally, we obtain the characteristic $2$ version of Theorem \ref{extensionvalues}.

\begin{theorem} \label{valueschar2} Let $n = 2m+1$ be odd and $q$ be even.  Let $y\tau$ be an element of $G_n^+ = G_n\langle \tau \rangle$ (or $U_n^+ = U_n\langle \tau \rangle$) such that $(y\tau)^2$ is regular unipotent, and let $\chi$ be a character of $G_n^+$ (or $U_n^+$) which is an extension of a real-valued irreducible character of $G_n$ (or $U_n$).  Then $\chi(y\tau) = \pm 1$ if $\chi(1)$ is odd, and $\chi(y\tau) = 0$ if $\chi(1)$ is even.  Also, $\chi(y\tau) \equiv \chi(\tau) \pmod {2}$ for any irreducible $\chi$ of $G_n^+$ (or $U_n^+$).
\end{theorem}
\begin{proof} Since $\sigma = w_0 \tau$, and $w_0 \in G_n$, we have $G_n^+ = G_n\langle \sigma \rangle$.  Since $w_0 \in \bar{G}_n^{\tilde{F}}$, we have $\bar{G}_n^{\tilde{F}}\langle \sigma \rangle = \bar{G}_n^{\tilde{F}} \langle \tau \rangle$, and since $\bar{G}_n^{\tilde{F}} \cong U_n$, we have $\bar{G}_n^{\tilde{F}}\langle \sigma \rangle \cong U_n^+$.  When $n$ is odd, the regular unipotent elements of $\bar{G}_n^F\langle \sigma \rangle$ (or $\bar{G}_n^{\tilde{F}}\langle \sigma \rangle$) correspond exactly to the elements of $G_n^+$ (or $U_n^+$) of the form $y\tau$ such that $(y\tau)^2 = u$ is regular unipotent in $G_n$ (or $U_n$).  If $\chi$ is an irreducible of $G_n^+$ (or $U_n^+$) which is extended from a real-valued irreducible character, then we know from Proposition \ref{regunipotentvalues} that $\chi(y\tau) = \pm 1$ or $0$.  From Theorem \ref{GLL}, we know that $\chi(u) = \pm 1$ when $\chi(1)$ is odd and $\chi(u) = 0$ when $\chi(1)$ is even.  Since we have
$$ \chi(y\tau)^2 \equiv \chi(u) \pmod {2\mathcal{O}},$$
and $\chi(y\tau) = \pm 1$ or $0$, then we must have $\chi(y\tau) = \pm 1$ when $\chi(1)$ is odd and $\chi(y\tau) = 0$ when $\chi(1)$ is even.  Since $\tau^2 = 1$, then $\chi(\tau)$ is an integer such that 
$$\chi(\tau)^2 \equiv \chi(1) \pmod {2}.$$
It follows that $\chi(y\tau) \equiv \chi(\tau) \pmod {2}$.
\end{proof}


\begin{thebibliography}{9}
\bibitem{alvis1}
D.~Alvis, The duality operation in the character ring of a finite Chevalley group, \emph{Bull. Amer. Math. Soc.} \textbf{1} (1979), 907--911.

\bibitem{alvis2}
D.~Alvis, Duality in the character ring of a finite Chevalley group, \emph{Proc. Symp. Pure Math.} (A.M.S.), \textbf{37} (1980), 353--357.

\bibitem{alvis3}
D.~Alvis, Duality and character values of finite groups of Lie type, \emph{J. Algebra}, \textbf{74} (1982), 211-222.

\bibitem{bantan}
E.~Bannai and H.~Tanaka, The decomposition of the permutation character $1^{GL(2n, q)}_{GL(n, q^2)}$, \emph{J. Algebra} \textbf{265} (2003), no. 2, 496--512.

\bibitem{Ca85}
R.~Carter, Finite groups of Lie type: conjugacy classes and complex characters.  John Wiley and Sons, 1985.

\bibitem{At}
J.H.~Conway, R.T.~Curtis, S.P.~Norton, R.A.~Parker, and R.A.~Wilson, Atlas of finite groups, Clarendon Press, Oxford, 1985.

\bibitem{curt}
C.W.~Curtis, Truncation and duality in the character ring of a finite group of Lie type, \emph{J. Algebra} \textbf{62} (1980), 320--332

\bibitem{dellusz}
P.~Deligne and G.~Lusztig, Representations of reductive groups over finite fields, \emph{Ann. of Math.} (2) \textbf{103} (1976), no. 1, 103--161.

\bibitem{digne}
F.~Digne, Descente de Shintani et restriction des scalaires, \emph{J. London Math. Soc. (2)} \textbf{59} (1999), no. 3, 867--880.

\bibitem{dignemichelu}
F.~Digne and J.~Michel, Foncteurs de Lusztig et caract\`eres des groupes lin\'eaires et unitaires sur un corps fini, \emph{J. Algebra} \textbf{107} (1987), no. 1, 217--255.

\bibitem{dmbook}
F.~Digne and J.~Michel, Representations of finite groups of Lie type. London Mathematical Society Student Texts, 21, Cambridge University Press, Cambridge, 1991.

\bibitem{dignemichel}
F.~Digne and J.~Michel, Groupes r\'eductifs non connexes, \emph{Ann. Sci. \'Ecole Norm. Sup. (4)} \textbf{27} (1994), no. 3, 345--406. 

\bibitem{ellers} E. W.~Ellers and W.~Nolte, Bireflectionality of orthogonal groups and symplectic groups,
\emph{Arch. Math. (Basel)} \textbf{39} (1982), 113-118.

\bibitem{enn1}
V.~Ennola, On the conjugacy classes of the finite unitary groups, \emph{Ann. Acad. Sci. Fenn. Ser. A I No.} \textbf{313} (1962), 13 pages.

\bibitem{enn2}
V.~Ennola, On the characters of the finite unitary groups, \emph{Ann. Acad. Sci. Fenn. Ser. A I No.} \textbf{323} (1963), 35 pages.

\bibitem{feit1}
W.~Feit, Extensions of cuspidal characters of ${\rm GL}_m(q)$, \emph{Publ. Math. Debrecen} \textbf{34} (1987), no. 3-4, 273--297.

\bibitem{fulman} J.~Fulman and R.~Guralnick, Conjugacy class properties of the extension
${\rm GL}(n,q)$ generated by the inverse transpose involution, \emph{J. Algebra} \textbf{275} (2004), 356-396.

\bibitem{gow3} 
R.~Gow,  On the Schur indices of characters of finite classical groups,
\emph{J. London Math. Soc.} \textbf{24} (1981), 135-147.

\bibitem{gow1}
R.~Gow, Properties of the characters of the finite general linear group related to the transpose-inverse involution, \emph{Proc. London Math. Soc. (3)} \textbf{47} (1983), no. 3, 493--506.

\bibitem{gow2}
R.~Gow, Two multiplicity-free permutation representations of the
general linear group $\GL(n,q^2)$ , \emph{Math. Z.}  \textbf{188} (1984), 45-54.
\bibitem{green}
J.A.~Green, The characters of the finite general linear groups, \emph{Trans. Amer. Math. Soc.} \textbf{80} (1955), no. 2, 402--447.

\bibitem{GLL}
J.A.~Green, G.I.~Lehrer, and G.~Lusztig, On the degrees of certain group characters, \emph{Quart. J. Math. Oxford} \textbf{27} (1976), 1--4.

\bibitem{hotspr}
R.~Hotta and T.A.~Springer, A specialization theorem for certain Weyl group representations and an application to the Green polynomials of unitary groups, \emph{Invent. Math.} \textbf{41} (1977), no. 2, 113--127.

\bibitem{kap} I. Kaplansky, Linear algebra and geometry. Allyn and Bacon, Boston, 1969.

\bibitem{Ka81}
N.~Kawanaka, Fourier transforms of nilpotently supported invariant functions on a finite simple Lie algebra, \emph{Proc. Japan Acad.} \textbf{57} (1981), 461--464.

\bibitem{kawan}
N.~Kawanaka, Generalized Gel'fand-Graev representations and Ennola duality, In \emph{Algebraic groups and related topics (Kyoto/Nagoya, 1983)}, 175--206, \emph{Adv. Stud. Pure Math.}, 6, North-Holland, Amsterdam, 1985.

\bibitem{luszsrin}
G.~Lusztig and B.~Srinivasan, The characters of the finite unitary groups, \emph{J. Algebra} \textbf{49} (1977), no. 1, 167--171.

\bibitem{Mac}
I.G.~Macdonald, Symmetric functions and Hall polynomials, Second Edition, With Contributions by A. Zelevinsky.  Oxford Mathematical Monographs, Oxford Science Publications, The Clarendon Press, Oxford University Press, New York, 1995.

\bibitem{malle}
G.~Malle, Generalized Deligne-Lusztig characters, \emph{J. Algebra} \textbf{159} (1993), no. 1, 64--97.

\bibitem{ohm}
Z.~Ohmori, On a Zelevinsky theorem and the Schur indices of the finite unitary groups, \emph{J. Math. Sci. Univ. Tokyo} \textbf{2} (1997), no. 2, 417--433.

\bibitem{shin}
T.~Shintani, Two remarks on irreducible characters of finite general linear groups, \emph{J. Math. Soc. Japan} \textbf{28} (1976), no. 2, 396--414.

\bibitem{sorlin1}
K.~Sorlin, Repr\'esentations de Gelfand-Graev pour les groupes r\'eductifs non connexes, \emph{C.R. Math. Acad. Sci. Paris} \textbf{334} (2002), no. 3, 179--184. 

\bibitem{sorlin2}
K.~Sorlin, \'El\'ements r\'eguliers et repr\'esentations de Gelfand-Graev pour les groupes r\'eductifs non connexes, \emph{Bull. Soc. Math. France} \textbf{132} (2004), no. 2, 157--199.

\bibitem{Spalt}
N.~Spaltenstein, Classes unipotentes et sous-groupes de Borel. Lecture Notes in Mathematics 946, Springer--Verlag, Berlin--New York, 1982.

\bibitem{St}
R.~Steinberg, Endomorphisms of linear algebraic groups,
\emph{Mem. Amer. Math. Soc.} {\bf 80} (1968), 1-108. 

\bibitem{TV}
N.~Thiem and C.R.~Vinroot, On the characteristic map of finite unitary groups, \emph{Adv. Math.} {\bf 210} (2007), no. 2, 707--732.

\bibitem{wall}
G.E.~Wall, On the conjugacy classes in the unitary, orthogonal, and symplectic groups, \emph{J. Austral. Math. Soc.} \textbf{3} (1962), 1--62.

\bibitem{won} M. J.~Wonenburger, 
 Transformations which are products of two involutions, \emph{J. Math. Mech.} \textbf{16} (1966), 113-118. 

\bibitem{zel}
A.V.~Zelevinsky, Representations of finite classical groups.  A Hopf algebra approach. Lecture Notes in Mathematics 869, Springer--Verlag, Berlin--New York, 1981.

\end{thebibliography}
\end{document}